\documentclass[final,3p,times]{elsarticle}

\biboptions{numbers,sort&compress}

\usepackage{hyperref}

\usepackage{amssymb}
\usepackage{amsmath}
\usepackage{amsthm}
\usepackage{soul}

\theoremstyle{remark}
\newtheorem{remark}{Remark}

\usepackage{comment}

\newcommand{\norm}[1]{\left\lVert#1\right\rVert}
\newcommand{\parens}[1]{\left( #1 \right)}
\newcommand{\curlies}[1]{\left\{ #1 \right\}}

\newcommand{\brackets}[1]{\left[ #1 \right]}

\newcommand{\Rbb}{\mathbb{R}}

\newcommand{\Bbf}{\mathbf{B}}
\newcommand{\cbf}{\mathbf{c}}

\newcommand{\gbf}{\mathbf{g}}
\newcommand{\Gbf}{\mathbf{G}}

\newcommand{\Ibf}{\mathbf{I}}

\newcommand{\kbf}{\mathbf{k}}

\newcommand{\nbf}{\mathbf{n}}

\newcommand{\rbf}{\mathbf{r}}

\newcommand{\Sbf}{\mathbf{S}}

\newcommand{\ubf}{\mathbf{u}}
\newcommand{\Ubf}{\mathbf{U}}
\newcommand{\vbf}{\mathbf{v}}
\newcommand{\Vbf}{\mathbf{V}}
\newcommand{\wbf}{\mathbf{w}}
\newcommand{\Wbf}{\mathbf{W}}
\newcommand{\xbf}{\mathbf{x}}

\newcommand{\zbf}{\mathbf{z}}
\newcommand{\Zbf}{\mathbf{Z}}

\newcommand{\rbft}{\mathbf{\tilde{r}}}
\newcommand{\ubft}{\mathbf{\tilde{u}}}

\newcommand{\mubf}{\boldsymbol{\mu}}
\newcommand{\muhat}{\hat{\boldsymbol{\mu}}}
\newcommand{\mustar}{\boldsymbol{\mu}^\star}

\newcommand{\Gen}{\mathcal{G}_{\text{en}}}
\newcommand{\Gde}{\mathcal{G}_{\text{de}}}

\newcommand{\Thetabf}{\mathbf{\Theta}}
\newcommand{\thetabf}{\boldsymbol{\theta}}

\newcommand{\psibf}{\boldsymbol{\psi}}

\newcommand{\Phibf}{\boldsymbol{\Phi}}
\newcommand{\dt}{\Delta t}

\newcommand{\phidot}{\dot{\phi}}

\newcommand{\Phidotbf}{\mathbf{\dot{\Phi}}}

\newcommand{\vecPhi}{\boldsymbol{\phi}}
\newcommand{\vecPhidot}{\boldsymbol{
\dot{\phi}}}
\newcommand{\ldbf}{\boldsymbol{\lambda}}
\newcommand{\omegabf}{\boldsymbol{\Omega}}

\usepackage{nicefrac}

\newcommand{\Dcal}{\mathcal{D}}
\newcommand{\Scal}{\mathcal{S}}

\newcommand{\Xcal}{\mathcal{X}}

\newcommand{\dd}[2]{\frac{d #1}{d #2}}

\newcommand{\pd}[2]{\frac{\partial #1}{\partial #2}}

\usepackage{xcolor}
\usepackage{tabularx}

\usepackage{booktabs}
\usepackage{multirow}
\usepackage{array}   % needed for >{\bfseries} column modifiers
\usepackage{siunitx}
\usepackage{xfp}     % provides \fpeval

\usepackage{subcaption}
\usepackage{caption}
\usepackage{threeparttable}

\begin{document}
\begin{frontmatter}


%\title{Accelerating Gradient-based Optimization for Time-dependent PDE Constraints via  Weak-form Latent Space Dynamics Identification}

%\title{Gradient-based Optimization for Time-dependent PDE Constraints using  Weak-form Reduced-order Modeling}

\title{Time-Dependent PDE-Constrained Optimization via 
\\ Weak-Form Latent Dynamics}


\author[APPMCU]{April Tran\corref{cor}}
\ead{chi.tran@colorado.edu}

\author[LLNL]{Terry Haut}
\ead{haut3@llnl.gov}

\author[APPMCU]{David M.~Bortz}
\ead{david.bortz@colorado.edu}

\author[LLNL]{Youngsoo Choi}
\ead{choi15@llnl.gov}

\affiliation[APPMCU]{organization={Department of Applied Mathematics, University of Colorado}, city={Boulder}, state={CO},  postcode={80309}, country={USA}}
\affiliation[LLNL]{organization={Center for Applied Scientific Computing, Lawrence Livermore National Laboratory}, city={Livermore}, postcode={94550}, state={CA},  country={USA}}

\cortext[cor]{Corresponding author}

%% Abstract
\begin{abstract}
Optimization problems constrained by high-dimensional, time-dependent partial differential equations require repeated forward and sensitivity solves, making high-fidelity optimization computationally prohibitive in many-query design and control settings. We present a weak-form latent-space reduced-order modeling framework for accelerating gradient-based PDE-constrained optimization. The proposed approach builds on Weak-form Latent Space Dynamics Identification (WLaSDI), which compresses high-dimensional solution trajectories into a low-dimensional latent representation and identifies parametric latent dynamics using weak-form system identification. By avoiding explicit numerical differentiation of training trajectories, the weak-form improves robustness to noisy data and yields more reliable surrogate dynamics for optimization. We formulate the resulting reduced PDE-constrained optimization problem and derive both direct-sensitivity and adjoint-based gradient expressions for the learned latent dynamics, enabling scalable gradient evaluation with respect to design parameters. The framework is demonstrated on three time-dependent benchmark problems: thermal radiative transfer for optimal hohlraum design, the two-stream instability Vlasov–Poisson system, and the inviscid Burgers equation. Across these examples, WLaSDI produces accurate optimal designs, remains robust under noisy training data, and delivers substantial computational savings, including speedups of up to five orders of magnitude relative to full-order optimization. These results demonstrate that weak-form latent dynamics provide an efficient and noise-robust surrogate foundation for gradient-based optimization of complex time-dependent PDE systems.

\begin{comment}
Solving large-scale optimization problems constrained by high-dimensional, time-dependent partial differential equations (PDEs) is computationally challenging, as it requires repeated solution of the governing equations.
The \emph{Latent Space Dynamics Identification} (LaSDI) framework, a data-driven, projection-based reduced order modeling method, addresses this challenge by enabling efficient surrogate-based simulation within optimization loops.
When paired with the weak formulation, the LaSDI framework gains improved accuracy and robustness to noise in training data.
In this work, we employ weak-form LaSDI to accelerate PDE-constrained optimization and derive adjoint-based sensitivities for scalable gradient evaluation.
We demonstrate the performance of the WLaSDI surrogates on several benchmark problems: 
thermal radiative transfer (TRT) for the optimal design of a hohlraum device, the two-stream instability Vlasov-Poisson system, and the inviscid Burgers equation. 
WLaSDI achieves robust and accurate performance under noisy training data, while delivering up to five orders of magnitude speedup over optimization with full-order models.
\end{comment}
\end{abstract}



%% Keywords
%\begin{keyword}



%\end{keyword}

\end{frontmatter}


\section{Introduction}
\label{sec:intro}
A practical digital twin for complex systems must support repeated solutions of design and control problems governed by high-fidelity physical models. In many applications, particularly in PDE-constrained optimization, this entails solving parametric, time-dependent partial differential equations (PDEs) within iterative optimization loops. The computational cost of repeatedly invoking high-fidelity PDE solvers, however, renders such many-query workflows prohibitively expensive. This fundamental challenge has motivated the development of surrogate-based optimization strategies aimed at reducing computational cost.

One common approach constructs surrogate models by interpolating selected quantities of interest (QoIs) \cite{QueipoHaftkaShyyEtAl2005ProgressinAerospaceSciences,ForresterKeane2009ProgressinAerospaceSciences}. While effective in some settings, these methods typically require retraining the response surface whenever the QoI changes. An alternative strategy is to accelerate the solution process itself through the use of \emph{reduced-order models} (ROMs), which approximate the underlying PDE dynamics at substantially lower computational cost.

In practice, \emph{projection-based reduced-order models} (pROMs) \cite{BennerGugercinWillcox2015SIAMRev} have proven particularly effective for accelerating the solution of high-dimensional parametric PDEs. 
%These methods have demonstrated success across a wide range of applications, including fluid dynamics \cite{CopelandCheungHuynhEtAl2022ComputerMethodsinAppliedMechanicsandEngineering, CheungChoiCopelandEtAl2023JournalofComputationalPhysics, MaulikLuschBalaprakash2020PhysFluids}, aerospace engineering \cite{AmsallemFarhat2008AIAAJournal,Bui-ThanhDamodaranWillcox2004AIAAJournal}, and design optimization \cite{ChoiBoncoraglioAndersonEtAl2020JournalofComputationalPhysics,AmsallemZahrChoiEtAl2015StructMultidiscOptim}, which motivates the focus of this study. 
%Within this class, non-intrusive pROMs are especially attractive, as they avoid intrusive modifications to legacy solvers and can be constructed using both simulation and experimental data.
In a typical pROM workflow, the high-dimensional full-order state is projected onto a low-dimensional subspace. This subspace is commonly constructed using techniques such as proper orthogonal decomposition (POD) \cite{BerkoozHolmesLumley1993AnnuRevFluidMech}, reduced basis methods \cite{RozzaHuynhPatera2008ArchComputatMethodsEng}, balanced truncation \cite{SafonovChiang1989IEEETransAutomatContr}, or nonlinear embeddings learned via autoencoders \cite{HintonSalakhutdinov2006Science}. To enable predictions at new parameter values, the resulting reduced model is typically interpolated in parameter space \cite{XieZhangWebster2019Mathematics,KadeethumBallarinChoiEtAl2022AdvancesinWaterResources}.
These methods have demonstrated success across a wide range of applications, including fluid dynamics \cite{CopelandCheungHuynhEtAl2022ComputerMethodsinAppliedMechanicsandEngineering, CheungChoiCopelandEtAl2023JournalofComputationalPhysics, MaulikLuschBalaprakash2020PhysFluids}, aerospace engineering \cite{AmsallemFarhat2008AIAAJournal,Bui-ThanhDamodaranWillcox2004AIAAJournal,McQuarrieHuangWillcox2021JournaloftheRoyalSocietyofNewZealand}, and design optimization \cite{ChoiBoncoraglioAndersonEtAl2020JournalofComputationalPhysics,AmsallemZahrChoiEtAl2015StructMultidiscOptim}, which motivates the focus of this study. 

Within the class of pROMs, the \emph{Latent Space Dynamics Identification (LaSDI)} framework \cite{FriesHeChoi2022ComputerMethodsinAppliedMechanicsandEngineering} has emerged as a strong contender due to its interpretability and non-intrusive nature.
Stemming from \cite{ChampionLuschKutzEtAl2019ProcNatlAcadSciUSA}, which employs Sparse Identification of Nonlinear Dynamics (SINDy) \cite{BruntonProctorKutz2016ProcNatlAcadSciUSA} to model latent space dynamics via an ordinary differential equation (ODE), LaSDI extends these ideas to a data-driven pROM framework for accelerating the simulation of parametric PDEs.
In essence, LaSDI projects the solution of a parametric PDE onto a low-dimensional parametric ODE governing the latent dynamics. By evolving this latent ODE in time and reconstructing the full-order state, LaSDI enables efficient and accurate simulation across the parameter domain at a substantially reduced computational cost. 

A key advantage of LaSDI lies in its interpretable and physics-informed representation of the latent dynamics. Rather than relying on a black-box neural network, the latent dynamics are modeled explicitly through a surrogate ODE, which promotes interpretability and facilitates the incorporation of physical constraints and structure. Parametric dependence is incorporated directly through the coefficients of the latent parametric ODE, rather than implicitly through learned latent trajectories, as seen in
\cite{XieZhangWebster2019Mathematics,KadeethumBallarinChoiEtAl2022AdvancesinWaterResources}.
%Moreover, LaSDI achieves enhanced robustness and reliability through \emph{weak-form} identification of the latent dynamics \cite{TranHeMessengerEtAl2024ComputerMethodsinAppliedMechanicsandEngineering,HeTranBortzEtAl2025NumericalMethEngineering}, which mitigates sensitivity to noise in training data and enables stable model discovery even in imperfect data regimes.
Since its introduction, LaSDI has evolved into a unified framework and has been applied across a wide range of problems.
Key developments include interpolation strategies 
%to enhance predictive capability across the parameter domain
\cite{FriesHeChoi2022ComputerMethodsinAppliedMechanicsandEngineering,HeChoiFriesEtAl2023JournalofComputationalPhysics,BonnevilleChoiGhoshEtAl2024ComputerMethodsinAppliedMechanicsandEngineering},
advances in autoencoder architectures 
%that better reveal underlying dynamics
\cite{ParkCheungChoiEtAl2024ComputMethodsApplMechEng,stephany2026higher,StephanyChoi2025,AndersonChungChoi2025,chung2026latent}, and improved robustness to noise in training data through the weak-form extension 
\cite{TranHeMessengerEtAl2024ComputerMethodsinAppliedMechanicsandEngineering,HeTranBortzEtAl2025NumericalMethEngineering}.  Section~\ref{sec:LaSDI} provides a brief overview of the LaSDI framework, while a comprehensive discussion is given in \cite{BonnevilleHeTranEtAl2024}.

Of particular relevance to the present work is the weak-form extension of LaSDI, WLaSDI \cite{TranHeMessengerEtAl2024ComputerMethodsinAppliedMechanicsandEngineering,HeTranBortzEtAl2025NumericalMethEngineering}.
Instead of identifying latent dynamics using SINDy \cite{BruntonProctorKutz2016ProcNatlAcadSciUSA}, which rely on approximating pointwise time derivatives and can therefore be highly sensitive to noise, WLaSDI employs the Weak-form Estimation of Nonlinear Dynamics (WENDy) \cite{BortzMessengerDukic2023BullMathBiol}.
This approach avoids explicit differentiation by projecting the dynamics onto a set of smooth, compactly supported test functions and applying integration by parts. As a result, WLaSDI can construct accurate and reliable surrogate models, which is particularly important in PDE-constrained optimization settings. Details of the WLaSDI algorithms used in this work are provided in Section~\ref{sec:WENDy}.

In this paper, we apply WLaSDI to optimization problems constrained by time-dependent PDEs. Specifically, this work makes three main contributions:
\begin{itemize}
    \item We formulate a WLaSDI-based reduced framework for optimization constrained by time-dependent PDEs, replacing repeated high-fidelity solves with efficient integration of learned latent dynamics.
    \item We derive direct-sensitivity and adjoint-based gradient expressions for the reduced latent dynamics, enabling scalable gradient evaluation as the parameter dimension increases \cite{YamaleevDiskinNielsen200812thAIAAISSMOMultidiscipAnalOptimConf,ChoiBoncoraglioAndersonEtAl2020JournalofComputationalPhysics,ZahrPersson2016JournalofComputationalPhysics}.
    \item We demonstrate, on Burgers’ equation, the Vlasov--Poisson system, and thermal radiative transfer, that WLaSDI enables accurate and efficient PDE-constrained optimization while exhibiting improved robustness to noisy training data compared with strong-form LaSDI and interpolation-based surrogates.
\end{itemize}

This paper is organized as follows. In Section~\ref{sec:problem_statement}, we present the problem statement and the fully discrete formulation of the PDE-constrained optimization problem. Section~\ref{sec:background} provides an overview of the WLaSDI framework. In Section~\ref{sec:PDE_CO}, we describe the application of WLaSDI to accelerate PDE-constrained optimization, derive the gradient and adjoint expressions. Numerical results are presented in Section~\ref{sec:results}.

\subsection{Problem Statement}
\label{sec:problem_statement}
Consider a parameterized, time-dependent dynamical system obtained from the spatial discretization of an underlying partial differential equation. The state of the system is represented by a finite-dimensional vector
$\ubf : [0,T] \times \Dcal \rightarrow \mathbb{R}^{N_\ubf},$ 
which evolves according to the semi-discrete system
\begin{equation}
\dd{\ubf}{t}(t, \mubf) = 
\boldsymbol{\mathcal{N}}\big(\ubf(t, \, \mubf), t, \, \mubf \big),
\qquad
\ubf(0,\mubf) = \gbf(\mubf),
\label{eq:spatial_discretized_PDE}
\end{equation}
where $\mubf \in \Dcal \subset \mathbb{R}^{N_\Dcal}$ denotes a finite-dimensional parameter vector, $\boldsymbol{\mathcal{N}}$ is the spatially discretized operator and $\gbf(\mubf)$ is the parameter-dependent initial condition. 

%Throughout this work, we assume that the parametric dependence enters exclusively through the initial condition. 

Temporal discretization of Eq.~\eqref{eq:spatial_discretized_PDE} using a general time-integration scheme yields a sequence of algebraic equations defining the discrete state evolution. Let $t_n = n\dt$, with $\dt = \frac{T}{N}$ and $n = 0,1,\dots,N$, and denote the discrete state by
$\ubf_n(\mubf) := \ubf(t_n, \mubf) \in \Rbb^{N_\ubf}.$
%$$\ubf_n(\mubf) := \ubf(t_n, \mubf) \in \Rbb^{N_\ubf}.$$
For notational clarity, we omit the explicit dependence of $\ubf_n$ on the parameter $\mubf$ in what follows. 
The fully discrete system can be expressed compactly in residual form as
\begin{equation}
    \begin{aligned}
    &\rbf_0 \big(\ubf_0, \ \mubf \big) = \ubf_0 - \gbf(\mubf)  \\
    &\rbf_n\big(\ubf_0,\ubf_1,\dots,\ubf_n, \ \mubf\big)  
    \end{aligned}
    \begin{aligned}
        & = \mathbf{0}_{N_\ubf},\\
        & = \mathbf{0}_{N_\ubf}, \qquad n = 1, 2, \dots, N.
    \end{aligned}
\label{eq:pde_residuals}
\end{equation}
where $\rbf_n: \Rbb^{N_\ubf} \times \Dcal \rightarrow \Rbb^{N_\ubf}$ denotes the discrete PDE residual at time step $n$. Note that $\rbf_n$ depends on the chosen spatial discretization and time-integration scheme and may involve the full state history, allowing for implicit or multistep methods.

The focus of this paper is the following fully discrete PDE-constrained optimization problem:
\begin{equation}
\boxed{
\begin{aligned}
\min_{\ubf_0,\ldots,\ubf_N, \ \mubf} \quad
& f(\ubf_0,\ldots,\ubf_N,\mubf) \\[0.5ex]
\text{s.t.} \quad
& \rbf_n(\ubf_0,\ldots,\ubf_n, \ \mubf) = \mathbf{0}_{N_\ubf},
\qquad n = 0,1,\ldots,N,\\
& \cbf(\ubf_0,\ldots,\ubf_N, \ \mubf) \le \mathbf{0}_{N_\cbf}. %\yc{\text{You mean } \mathbf{0}_{N_\cbf} ?}.
\end{aligned}
}
\label{eq:PDE_CO_discrete}
\end{equation}
Here, $f$ is a scalar-valued objective functional that depends on the discrete state trajectory and the parameter vector $\mubf$, and $\cbf \in \mathbb{R}^{N_\cbf}$ represents a collection of inequality constraints. The optimization seeks parameters $\mubf$ (and the corresponding state trajectory $\{\ubf_n\}_{n=0}^N$ satisfying the discrete PDE constraints in Eq.~\eqref{eq:pde_residuals}) that minimize the objective while respecting the inequality constraints.

\section{Background}
\label{sec:background}
\subsection{LaSDI: Latent Space Dynamics Identification}
\label{sec:LaSDI}
In this section, we provide an overview of the Latent Space Dynamics Identification (LaSDI) framework \cite{FriesHeChoi2022ComputerMethodsinAppliedMechanicsandEngineering}. LaSDI is a data-driven simulation methodology designed to accelerate the numerical solution of high-dimensional, parameterized partial differential equations (PDEs). The central idea of LaSDI is to project the governing parametric, semi-discrete system in Eq.~\eqref{eq:spatial_discretized_PDE} onto a low-dimensional latent space and to identify a parametric ODE that governs the evolution of the latent variables. In doing so, LaSDI effectively recasts the original problem from repeatedly solving a high-dimensional PDE to integrating a low-dimensional ODE. 

We briefly recall the parametric system in Eq.~\eqref{eq:spatial_discretized_PDE} and its fully discrete numerical solution $\ubf_n(\mubf)$ for $n= 0, 1, \dots, N$. The discrete trajectory is collected as: 
$$\Ubf(\mubf) := \left[\begin{array}{ccc}
         \ubf_0(\mubf) &  \cdots & \ubf_N(\mubf)
         \end{array}\right]^T \in \Rbb^{(N+1) \times N_\ubf}. $$

\begin{figure}[htbp]
    \centering
    \includegraphics[width=0.9\textwidth]{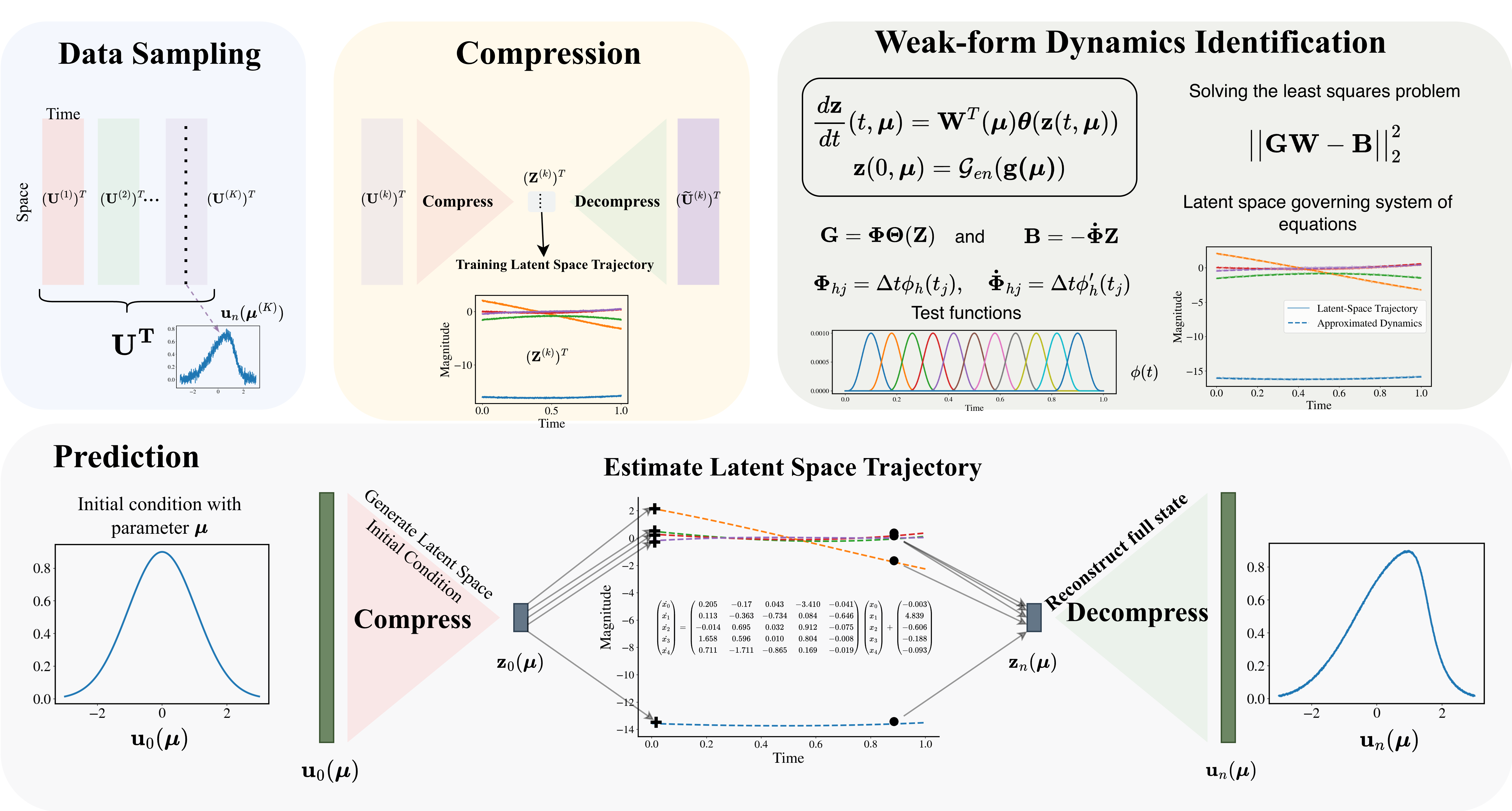}
    \caption{Overview of the WLaSDI framework, featuring 4 key stages: Data Sampling, Compression, Weak-form Dynamics Identification, and Prediction }
    \label{fig:WLaSDI}
\end{figure}

We now describe the LaSDI framework in its most general form, which consists of four key stages: data sampling, compression, latent space dynamics identification, and prediction. 
A comprehensive overview of the method can be found in \cite{BonnevilleHeTranEtAl2024}. A schematic illustrating the workflow of the weak form extension, WLaSDI, is provided in Figure~\ref{fig:WLaSDI}.

\paragraph{\textbf{Data Sampling}}
In this stage, a finite set of training parameters is selected from the parameter domain $\Dcal$. For each selected parameter value, a high-fidelity numerical solution of Eq.~\eqref{eq:spatial_discretized_PDE} is computed to generate training trajectories. The selection of training parameters can follow either uniform sampling or an adaptive procedure; see \cite{BonnevilleHeTranEtAl2024,HeChoiFriesEtAl2023JournalofComputationalPhysics,HeTranBortzEtAl2025NumericalMethEngineering,ChoiBoncoraglioAndersonEtAl2020JournalofComputationalPhysics} for details.

Let $\mathcal{S} = \curlies{ \mubf^{(1)}, \mubf^{(2)}, \dots, \mubf^{(K)} }  \subset \mathcal{D}$  denote the set of training parameters. For each $\mubf^{(k)} \in \mathcal{S}$, the computed trajectories is denoted by $\Ubf(\mubf^{(k)}):= \Ubf^{(k)} \in \Rbb^{(N+1) \times N_\ubf}$. All training trajectories are then concatenated to form the global training data matrix 
$\Ubf := \left[\begin{array}{ccc}
         \parens{\Ubf^{(1)}}^T &  \cdots & \parens{\Ubf^{(K)}}^T
         \end{array}\right]^T \in \Rbb^{(N+1)K \times N_\ubf}. $

\paragraph{\textbf{Compression}} 
Given the collection of high-dimensional training trajectories, the next stage of LaSDI seeks to obtain a low-dimensional representation of the system state. 
Let $\zbf(t, \mubf) \in \Rbb^{N_\zbf}$ denote the latent variable associated with 
$\ubf(t, \mubf)$, where $N_\zbf \ll N_\ubf$ is the latent dimension. The mapping between the full-order and latent representations is defined through an encoder–decoder pair,
\begin{equation*}
    \begin{aligned}
        \Gen: \Rbb^{N_{\ubf}}  \rightarrow \Rbb^{N_{\zbf}} ,
        \qquad
 \Gde: \Rbb^{N_{\zbf}}   \rightarrow \Rbb^{N_{\ubf}}  .
    \end{aligned}
\end{equation*}
These operators may be constructed using Proper Orthogonal Decomposition (POD) or autoencoder neural networks. For notational convenience, we define 
$\zbf_n(\mubf) := \Gen(\ubf_n(\mubf))$, and conversely, 
$\ubf_n(\mubf) = \Gde(\zbf_n(\mubf)).$
The full-time-discretized latent space trajectory 
%$\Zbf(\mubf) \in \Rbb^{(N+1) \times N_\zbf} $
is then given by 
$$\Zbf(\mubf) := \left[\begin{array}{ccc}
         \zbf_0(\mubf) &  \cdots & \zbf_N(\mubf)
         \end{array}\right]^T \in \Rbb^{(N+1) \times N_\zbf}. $$

\paragraph{\textbf{Dynamics Identification}}
\label{par:DI}
Within the LaSDI framework, the temporal evolution of the latent variable $\zbf(t, \mubf) \in \Rbb^{N_\zbf}$ for an arbitrary $\mubf$ is assumed to satisfy a parametric surrogate ODE of the form
\begin{equation}
    \dd{\zbf}{t} \ (t, \mubf) = \Wbf^T(\mubf) \ \thetabf \parens{\zbf(t, \mubf) }, 
    \qquad \zbf(0, \mubf) = \Gen(\gbf(\mubf)),
     \label{eq:paramtric_latent_ode}
\end{equation}
where $\thetabf(\zbf)  :=  \left[\begin{array}{ccc}
f_1(\zbf) &  \cdots &  f_J(\zbf) \end{array}\right]^T \in \mathbb{R}^{{J}}$ denotes the user-specified library of candidate nonlinear functions $f_j: \Rbb^{N_\zbf} \rightarrow \Rbb$. The coefficient matrix $\Wbf(\mubf) \in \Rbb^{J \times N_\zbf}$ encodes the parametric dependence of the latent dynamics.

In the offline phase of LaSDI, latent trajectories corresponding to the training parameters 
%$ \Scal = \curlies{\mubf^{(k)}}_{k=1}^K$
$\Scal$
are used to identify the coefficient matrices governing the latent dynamics. Specifically, for each $\mubf^{(k)} \in \Scal$, the corresponding coefficient matrix $\Wbf^{(k)} := \Wbf(\mubf^{(k)})$  is identified from the latent trajectory $\Zbf^{(k)} := \Zbf(\mubf^{(k)})$; see section \ref{sec:WENDy} for details\footnote{Rather than performing compression and dynamics identification sequentially, as described here, the autoencoder and latent space ODE can alternatively be trained simultaneously in an end-to-end manner; see \cite{HeChoiFriesEtAl2023JournalofComputationalPhysics, BonnevilleChoiGhoshEtAl2024ComputerMethodsinAppliedMechanicsandEngineering} for details.}.
In the online phase, \hyperref[par:prediction]{Prediction}, given $\mubf \in \Dcal$, the associated coefficient matrix $\Wbf(\mubf)$ is obtained by interpolating the precomputed training coefficients. 

%$\curlies{ \Wbf^{(k)}}^K_{k = 1}.$A range of interpolation strategies is discussed in \ref{sec:Interpolation}.

This formulation offers several advantages. First, it is highly flexible: the choice of the library $\thetabf$ is user-defined and can be tailored to incorporate structural or physical constraints in the latent space. Second, by interpolating coefficient matrices rather than latent trajectories, parametric generalization is achieved efficiently and robustly. Finally, the resulting model is interpretable, as the latent dynamics are represented explicitly as a parametric ODE rather than a black-box surrogate.

%using SINDy \cite{BruntonProctorKutz2016ProcNatlAcadSciUSA} or  WENDy \cite{BortzMessengerDukic2023BullMathBiol}; see section \ref{sec:WENDy} for details\footnote{Rather than performing compression and dynamics identification sequentially, as described here, the autoencoder and latent space ODE can alternatively be trained simultaneously in an end-to-end manner; see \cite{HeChoiFriesEtAl2023JournalofComputationalPhysics, BonnevilleChoiGhoshEtAl2024ComputerMethodsinAppliedMechanicsandEngineering} for details.}.

\paragraph{\textbf{Prediction}}
\label{par:prediction}
Given an arbitrary parameter value $\mubf \in \Dcal$, the corresponding coefficient matrix $\Wbf(\mubf)$ is obtained via interpolation of the precomputed coefficients $\curlies{ \Wbf^{(k)}}^K_{k = 1}.$ A range of interpolation strategies for this purpose is discussed in \ref{sec:Interpolation}.

Following Eq.~\eqref{eq:paramtric_latent_ode}, the latent initial condition is obtained by encoding the full-order initial condition,
$\zbf_0(\mubf) = \zbf(0, \mubf) = \Gen(g(\mubf)).$
%$$\zbf_0(\mubf) = \zbf(0, \mubf) = \Gen(g(\mubf)).$$
With $\Wbf(\mubf)$ and $\zbf_0(\mubf)$ determined, the full latent space trajectory $\Zbf(\mubf)$ is obtained by numerically solving the identified IVP. 
Finally, the latent trajectory is decoded to recover the full-order solution $\widetilde{\Ubf}(\mubf)$. %$\ubf_n = \Gde(\zbf_n, \mubf)$. 

\begin{comment}
Let $\omegabf$ denote a one-step\footnote{note that for} time integration scheme, which may be either implicit or explicit. Then, the time-discrete evolution of the latent variable satisfies
$$\psibf \parens{\zbf_n(\mubf),  \zbf_{n-1}(\mubf), \mubf} = \mathbf{0} \in \Rbb^{N_\zbf}, \text{ for } n = 1, \dots, N;$$
\end{comment}

\subsection{Weak-form Identification of Latent Space Dynamics}
\label{sec:WENDy}
This section focuses on the identification of latent space dynamics within the LaSDI framework, specifically the coefficient matrix $\Wbf$ appearing in the surrogate latent ODE in the \hyperref[par:DI]{Dynamics Identification} stage, using the weak-form. This approach is rooted in the \emph{weak-form} of a differential equation, obtained by multiplying the governing equation by a smooth test function $\phi$ and applying integration by parts to transfer derivatives from the state to $\phi$, yielding an equivalent expression involving fewer derivatives. In the context of system identification, the weak-form therefore avoids the explicit approximation of time derivatives from potentially noisy data.

Weak-form approaches to system identification have a long history, dating back to early studies in the 1950s \cite{Shinbrot1957JFluidsEng} and 1960s \cite{LoebCahen1965IEEETransAutomatContr}.  For a review of the advances up until the early 1990s, see \cite{PreisigRippin1993ComputersChemicalEngineering}. More recently, following the introduction of the Sparse Identification of Nonlinear Dynamics (SINDy) \cite{BruntonProctorKutz2016ProcNatlAcadSciUSA, RudyBruntonProctorEtAl2017SciAdv} framework, weak-form variants of SINDy have been shown to substantially improve robustness to noise by avoiding unstable derivative approximations \cite{PantazisTsamardinos2019Bioinformatics, GurevichReinboldGrigoriev2019ChaosInterdiscipJNonlinearSci, SchaefferMcCalla2017PhysRevE, WangHuanGarikipati2019ComputerMethodsinAppliedMechanicsandEngineering, ReinboldGurevichGrigoriev2020PhysRevE}.
More broadly, one can consider \emph{Weak-form Scientific Machine Learning} (WSciML) framework \cite{MessengerTranDukicEtAl2024SIAMNews} as encompassing both equation discovery (WSINDy) and parameter estimation (WENDy, see below). WSciML builds on methods from sparse regression, total least squares, and maximum likelihood, and is supported by theory on asymptotic consistency \cite{MessengerBortz2024IMAJNumerAnal} and optimal test-function design \cite{TranBortz2026SISC, MessengerBortz2021MultiscaleModelSimul,MessengerBortz2021JournalofComputationalPhysics,RummelMessengerBeckerEtAl2025,BortzMessengerDukic2023BullMathBiol,ChawlaBortzDukic2026HandbookofVisualExperimentalandComputationalMathematics}. Notably, the use of WSciML has been successfully applied to real world problems in biology \cite{MessengerWheelerLiuEtAl2022JRSocInterface,MessengerDwyerDukic2024JRSocInterface,Heitzman-BreenDukicBortz2026BullMathBiol,LyonsDukicBortz2025PLoSComputBiol,MinorElderdVanAllenEtAl2025}, atmospheric science \cite{MinorMessengerDukicEtAl2025JournalofGeophysicalResearchMachineLearningandComputation}, and plasma physics 
\cite{VaseyMessengerBortzEtAl2025JournalofComputationalPhysics,MessengerSouthworthHammerEtAl2025}.
Within reduced-order modeling, these ideas form the basis for the Weak-form LaSDI (WLaSDI) algorithms \cite{TranHeMessengerEtAl2024ComputerMethodsinAppliedMechanicsandEngineering, HeTranBortzEtAl2025NumericalMethEngineering}, which apply a WSciML methodology to latent space dynamics, enabling efficient and robust surrogate modeling. 

We now present the weak-form identification of the latent space ODE that forms the basis of the WLaSDI algorithm. To illustrate, we reconsider latent space dynamics, omitting parametric dependence in Eq.~\eqref{eq:paramtric_latent_ode}. The latent state $\zbf(t) \in \Rbb^{N_\zbf}$ is assumed to satisfy the ODE
\begin{equation}
        \dd{\zbf}{t} (t) = \Wbf^T \thetabf \parens{\zbf(t) }, 
    \qquad \zbf(0) = \zbf_0,
    \label{eq:latent_ode}
\end{equation}
where $\thetabf(\zbf) \in \Rbb^J$ denotes a library of candidate functions, and $\Wbf \in \Rbb^{J \times N_\zbf}$ is the unknown coefficient matrix. 
Suppose the latent trajectory $\curlies{\zbf_n}_{n=0}^N$, with $\zbf_n = \zbf(t_n)$, are available at uniformly spaced time instances $t_n = n\dt$. 
In the original LaSDI implementation, $\Wbf$ is identified using SINDy \cite{BruntonProctorKutz2016ProcNatlAcadSciUSA}, which solves the linear system\footnote{While SINDy is typically formulated as a sparse regression problem, sparsity is not enforced in the LaSDI setting, as the goal is to construct a surrogate model of the latent dynamics.}
$$\dot{\Zbf} \approx  \Thetabf(\Zbf) \Wbf,$$
where $\dot{\Zbf} \in \Rbb^{(N+1) \times N_\zbf}$ and $\Thetabf(\Zbf) \in \Rbb^{(N+1) \times J}$ are defined as
\begin{equation*}
\dot{\Zbf}  :=  \left[\begin{array}{c}
    \dot{\zbf}_0^T \\  \vdots \\ \dot{\zbf}_N^T
         \end{array}\right], 
    \qquad
 \Thetabf(\Zbf) := \left[\begin{array}{ccc}
f_1(\zbf_0) &  \cdots &  f_J(\zbf_0)\\
\vdots &  \ddots &  \vdots \\
f_1(\zbf_N) &  \cdots &  f_J(\zbf_N)\\
\end{array}\right].
\end{equation*}

Computing $\dot{\Zbf}$ directly requires numerical differentiation of the latent trajectories, typically via finite difference approximations, which can be inaccurate and highly sensitive to noise. To avoid explicit differentiation, WLaSDI instead adopts the weak-form for identifying the coefficient matrix $\Wbf$. Multiplying both sides of Eq.~\eqref{eq:latent_ode} by a smooth, compactly supported test function $\phi(t)$ and integrating by parts yields
 \begin{equation}
    - \int_0^T \dot{\phi}  \zbf dt = \int_0^T \phi \Wbf^T \thetabf(\zbf) dt.
     \label{eq:weak_ode}
\end{equation}
Approximating these integrals using the trapezoidal rule leads to
$ - \boldsymbol{\dot{{\phi}}}^T \Zbf \approx  \boldsymbol{\phi}^T\mathbf{\Theta(Z)} \mathbf{W},
$    
where the row vectors are defined as
\begin{equation*}
\begin{aligned}
    \boldsymbol{\phi}^T  := \left[\begin{array}{ccc}
         \phi(t_0) &  \cdots &  \phi(t_N)
         \end{array}\right]\mathcal{Q}, 
         \qquad
     \boldsymbol{\dot{{\phi}}}^T := \left[\begin{array}{ccc}
         \phidot(t_0) &  \cdots &  \phidot(t_N)
         \end{array}\right]\mathcal{Q}, 
    \end{aligned}
    \label{eq:vecPhi}
\end{equation*}
with quadrature matrix $\mathcal{Q} = \text{diag}\parens{\frac{\dt}{2}, \dt, \cdots, \dt, \frac{\dt}{2}} \in \Rbb^{(N+1) \times (N+1)}.$ Using  a collection of $M$ test functions, $\mathcal{M} := \curlies{\phi_m}_{m=1}^M$, the identification of $\mathbf{W}$ can be cast as the weak-form least squares problem:
\begin{equation}
    \min_{\mathbf{W}} \norm{ \mathbf{B} -\mathbf{G}\mathbf{W} }_2^2,
    \label{eq:WENDy}
\end{equation}
where the matrices 
$\Gbf \in  \Rbb^{M \times J},$
$\Bbf \in \Rbb^{M \times N_\zbf},$
along with $\Phibf, \Phidotbf \in \Rbb^{M \times (N+1)}$, are defined by
\begin{equation*}
\begin{aligned}
\Gbf := \Phibf \Thetabf(\Zbf), 
\quad 
\Bbf := -\Phidotbf \Zbf, 
\qquad\textrm{with}\qquad
\Phibf :=  \left[\begin{array}{ccc} 
\vecPhi_1^T & \cdots & \vecPhi_M^T
 \end{array}\right]^T, 
 \quad
 \Phidotbf & := \left[\begin{array}{ccc} 
\vecPhidot_1^T & \cdots & \vecPhidot_M^T
 \end{array}\right]^T.
\end{aligned}
\end{equation*}
This is formally referred to as \emph{Weak-form Estimation of Nonlinear Dynamics (WENDy)} \cite{BortzMessengerDukic2023BullMathBiol} within the WSciML framework\footnote{WENDy focuses on parameter estimation, in contrast to Weak form Sparse Identification of Nonlinear Dynamics (WSINDy) \cite{MessengerBortz2021JournalofComputationalPhysics, MessengerBortz2021MultiscaleModelSimul} that target \emph{sparse} equation discovery.}.
The coefficient matrix $\Wbf$ can be computed using ordinary least squares (WENDy-OLS), or iteratively reweighted least squares (WENDy-IRLS) to account for the errors-in-variables nature of Eq.~\eqref{eq:WENDy}, or maximum-likelihood estimation (WENDy-MLE) for nonlinear-in-parameter models \cite{RummelMessengerBeckerEtAl2025}. 
The choice of test function collection $\mathcal{M}$ plays a critical role in weak-form inference; strategies for test function design are discussed in \cite{TranBortz2026SISC, BortzMessengerDukic2023BullMathBiol}.

By using WENDy to obtain the governing equation of the latent space, the \emph{Weak form Latent Space Dynamics Identification (WLaSDI)} algorithm \cite{TranHeMessengerEtAl2024ComputerMethodsinAppliedMechanicsandEngineering,HeTranBortzEtAl2025NumericalMethEngineering} significantly improves the robustness of LaSDI to noise in the training data.
Figure~\ref{fig:WLaSDI} presents a schematic overview of the WLaSDI workflow.
WLaSDI has been shown to enable fast, robust, and accurate surrogate simulations, achieving up to three orders of magnitude noise reduction compared to the strong form while providing speedups of two orders of magnitude over high-fidelity solvers. These improvements yield more reliable and efficient sensitivity and gradient evaluations, directly accelerating convergence in gradient-based PDE-constrained optimization, which motivates our use of WLaSDI in this work.

\section{Accelerating Gradient-based PDE-Constrained Optimization with WLaSDI Framework }
\label{sec:PDE_CO}

\subsection{Gradient-based Optimization with Time-Dependent PDE Constraints}
\label{sec:gradient_optimzation}
This section focuses on computing the gradients required for gradient-based optimization of the time-dependent PDE-constrained optimization problem defined in Eq.~\eqref{eq:PDE_CO_discrete}. For clarity of notation, we omit the explicit parametric dependence of the state vectors and write $\ubf_n \equiv \ubf_n(\mubf)$.

Gradient-based optimization algorithms require evaluation of the first derivatives of the objective functional
$f(\ubf_0,\dots,\ubf_N, \ \mubf)$ and any constraint functions
$\cbf(\ubf_0,\dots,\ubf_N, \ \mubf)$ with respect to the parameter vector
$\mubf \in \Dcal \subset \mathbb{R}^{N_\Dcal}$. Without loss of generality, we consider the derivative of the objective function $f$ with respect to an individual parameter component $\mu_i$, $\dd{f}{\mu_i}$, for $i = 1,\dots,N_\Dcal$. Applying the chain rule yields
\begin{equation}
    \begin{aligned}
    \dd{f}{\mu_i}(\ubf_0, \dots, \ubf_N, \mubf) 
    =  \sum_{n=0}^N \pd{f}{\ubf_n} (\ubf_0, \dots, \ubf_N, \ \mubf) \  \pd{\ubf_n}{\mu_i}(\mubf) 
    + \pd{f}{\mu_i} (\ubf_0, \dots, \ubf_N,  \  \mubf). 
    \end{aligned}
    \label{eq:gradient_f}
\end{equation}
For a given objective function $f$, the derivatives $\pd{f}{\mu_i} \in \Rbb$ and 
$\pd{f}{\ubf_n} \in \Rbb^{1 \times N_\ubf}, n = 0, 1, \dots, N$  are generally straightforward.  The main challenge lies in evaluating the state sensitivities $\pd{\ubf_n}{\mu_i} \in \Rbb^{N_\ubf}$. To this end, we note that the discretized PDE residuals
$\rbf_n(\ubf_0, \dots, \ubf_n, \ \mubf)$, defined in Eq.~\eqref{eq:pde_residuals}, satisfy
%$\dd{\rbf_n}{\mu_i} = \mathbf{0}, n = 0, 1, \dots, N.$ Applying the chain rule gives
\begin{equation}
    \begin{aligned}
    %\dd{\rbf_0}{\mu_i}  &= \pd{\rbf_0}{\mu_i} + \pd{\rbf_0}{\ubf_0}\pd{\ubf_0}{\mu_i} = \mathbf{0}, \\
    \dd{\rbf_n}{\mu_i}  &= \pd{\rbf_n}{\mu_i}
    +
    \sum_{j=0}^n \pd{\rbf_n}{\ubf_{j}}  \ \pd{\ubf_j}{\mu_i}   = \mathbf{0}_{N_\ubf}, \qquad    n = 0, 1, \dots, N,
    \end{aligned}
    \label{eq:derivative_rn}
\end{equation}
where $\pd{\rbf_n}{\ubf_{j}} \in \Rbb^{N_\ubf \times N_\ubf}. $
Solving these relations forward in time yields the state sensitivities
\begin{equation}
    \begin{aligned}
        \pd{\ubf_0}{\mu_i} &= - \brackets{\pd{\rbf_0}{\ubf_0}}^{-1}\pd{\rbf_0}{\mu_i}, \\
        \pd{\ubf_{n}}{\mu_i} &=  - \brackets{\pd{\rbf_n}{\ubf_{n}}}^{-1}
        \parens{ \pd{\rbf_n}{\mu_i} +
        \sum_{j=0}^{n-1} \pd{\rbf_n}{\ubf_{j}} \ \pd{\ubf_{j}}{\mu_i}}, 
        \qquad n = 1, 2, \dots, N. 
    \end{aligned}
    \label{eq:direct}
\end{equation}
This procedure is referred to as the \emph{direct} sensitivity approach, and the resulting sensitivities can be substituted into the gradient expression in Eq.~\eqref{eq:gradient_f}.
Alternatively, the \emph{adjoint} approach can be derived by introducing Lagrange multipliers $\ldbf_n \in \Rbb^{N_\ubf}$. Using Eq.~\eqref{eq:derivative_rn}, we write  
\begin{equation}
    \begin{aligned}
        \sum_{n=0}^N \pd{f}{\ubf_n}  \
        \pd{\ubf_n}{\mu_i}
        %& =  \sum_{n=0}^N \pd{f}{\ubf_n}  \pd{\ubf_n}{\mu_i}  -  \sum_{n=0}^N \ldbf_n^T\underbrace{\dd{\rbf_n}{\mu_i}}_{\mathbf{0}}\\
        & = 
        \sum_{n=0}^N \pd{f}{\ubf_n}  \ 
        \pd{\ubf_n}{\mu_i} 
        -    \sum_{n=0}^N \ldbf_n^T \parens{ 
        \pd{\rbf_n}{\mu_i} +
    \sum_{j=0}^n \pd{\rbf_n}{\ubf_{j}} \ \pd{\ubf_j}{\mu_i}}.
    \end{aligned}
    \label{eq:adjoint_derivation}
\end{equation}
Rearranging terms yields
\begin{equation*}
\begin{aligned}
    \sum_{n=0}^N \pd{f}{\ubf_n} \ \pd{\ubf_n}{\mu_i} = 
    -\sum_{n=0}^N \ldbf_n^T \ \pd{\rbf_n}{\mu_i}
    + \sum_{n=0}^N \parens{
    \pd{f}{\ubf_n} - \sum_{j=n}^N \ldbf_j^T \ \pd{\rbf_j}{\ubf_n}
    }\pd{\ubf_n}{\mu_i}.
\end{aligned}
\end{equation*}
By enforcing the adjoint equations such that the expressions in parentheses vanish, $ \pd{f}{\ubf_n} - \sum_{j=n}^N \ldbf_j^T \ \pd{\rbf_j}{\ubf_n}  = \mathbf{0}^T_{N_\ubf},$
%\begin{equation*}begin{aligned}
    %\pd{f}{\ubf_n} - \sum_{j=n}^N \ldbf_j^T \ \pd{\rbf_j}{\ubf_n}  = \mathbf{0}^T_{N_\ubf}, \qquad n = 0, 1, \dots, N.
%\end{aligned}\end{equation*}
we can solve for the Lagrange multipliers $\ldbf_n \in \Rbb^{N_\ubf}$ backward in time by 
\begin{equation}
    \begin{aligned}
    \ldbf_N^T & = \pd{f}{\ubf_N} \brackets{\pd{\rbf_N}{\ubf_N}}^{-1},\\
        \ldbf_n^T  &= \parens{\pd{f}{\ubf_n} - \sum_{j=n+1}^N \ldbf_j^T \ \pd{\rbf_j}{\ubf_n}} 
        \brackets{\pd{\rbf_n}{\ubf_n}}^{-1}, 
        \qquad n = (N-1), (N-2), \dots, 0. 
    \end{aligned}
    \label{eq:adjoint}
\end{equation}
The gradient of the objective function $f$ in Eq.~\eqref{eq:gradient_f} then becomes 
\begin{equation}
\begin{aligned}
   \dd{f}{\mu_i}(\ubf_0, \dots, \ubf_N, \mubf)  =  \pd{f}{\mu_i}(\ubf_0, \dots, \ubf_N, \ \mubf) 
   - \sum_{n=0}^N \ldbf^T_n \ \pd{\rbf_n}{\mu_i}(\ubf_0, \dots, \ubf_N, \ \mubf) .
\end{aligned}
\label{eq:adjoint_gradient}
\end{equation}

\begin{remark}
The direct sensitivity approach requires solving a linear system at each time step for every parameter component $\mu_i$, resulting in a total of $(N+1)N_\Dcal$ linear solves (Eq.~\eqref{eq:direct}). Once these state sensitivities are computed, no additional solves are needed to evaluate gradients of the constraint functions $\cbf \in \mathbb{R}^{N_\cbf}$.
In contrast, the adjoint approach eliminates the dependence on the parameter dimension by introducing the Lagrange multiplier $\ldbf_n$ associated with the PDE constraints at each time step. This requires solving $(N+1)$ adjoint systems for the objective function and an additional $(N+1)$ systems per constraint, for a total of $(N+1)(1+N_\cbf)$ linear solves (Eq.~\eqref{eq:adjoint}). Consequently, the adjoint formulation is computationally advantageous when $1 + N_\cbf < N_\Dcal$.
\end{remark}

Using the gradients derived above, the PDE-constrained optimization problem in Eq.~\eqref{eq:PDE_CO_discrete} can be solved using standard nonlinear optimization algorithms, such as trust-region methods \cite{ConnGouldToint2000}, interior-point methods \cite{WachterBiegler2006MathProgram}, or Sequential Quadratic Programming (SQP) \cite{GillMurraySaunders2005SIAMRev}.  
However, the treatment of the discretized PDE constraints is computationally intensive. In particular, evaluating the state sensitivities $\pd{\ubf_n}{\mu_i} \in \mathbb{R}^{N_\ubf}$ at each time step requires solving large-scale linear systems associated with the full-order model, resulting in significant computational cost. To address this challenge, we employ Weak-form Latent Space Dynamics Identification (WLaSDI) to accelerate the optimization procedure by enabling efficient surrogate simulation of the underlying PDE. The resulting reduced optimization framework is presented in Section~\ref{sec:WLaSDI_optimization}.

\subsection{WLaSDI for PDE-Constrained Optimization}
\label{sec:WLaSDI_optimization}
We begin by recalling the full-order model (FOM) describing the evolution of the state
$\ubf(t,\mubf)$, governed by the semi-discrete system in
Eq.~\eqref{eq:spatial_discretized_PDE}.
Within the LaSDI framework (see Section~\ref{sec:background}), the high-dimensional
state $\ubf \in \mathbb{R}^{N_\ubf}$ is approximated by a low-dimensional latent
representation $\zbf \in \mathbb{R}^{N_\zbf}$, with $N_\zbf \ll N_\ubf$, through
encoder--decoder maps $\Gen$ and $\Gde$.
The evolution of the latent state $\zbf(t,\mubf)$ is governed by the learned surrogate
dynamics defined in Eq.~\eqref{eq:paramtric_latent_ode}.   

Applying a time-integration scheme to the latent dynamical system yields a sequence of \emph{reduced residuals}, denoted by $\rbft_n$, which enforce the discrete-time evolution of the latent variables and serve as reduced-order analogues of the full-order residuals
$\rbf_n$ defined in Eq.~\eqref{eq:pde_residuals}.
For notational convenience, we write $\zbf_n \equiv \zbf_n(\mubf)$.
The reduced residuals are given by
\begin{equation}
    \begin{aligned}
    &\rbft_0 \big(\zbf_0, \ \mubf \big) = \zbf_0 - \Gen(\gbf(\mubf))  \\
    &\rbft_n \big(\zbf_0,\zbf_1,\dots,\zbf_n, \ \mubf\big)  
    \end{aligned}
    \begin{aligned}
        & = \mathbf{0}_{N_\zbf},\\
        & = \mathbf{0}_{N_\zbf}, \qquad n = 1, 2, \dots, N.
    \end{aligned}
\label{eq:reduced_residuals}
\end{equation}
Consequently, the optimization in Eq.~\eqref{eq:PDE_CO_discrete} is reformulated as 
\begin{equation}
\boxed{
\begin{aligned}
\min_{\zbf_0,\dots,\zbf_N, \ \mubf} \quad
& f( \Gde(\zbf_0),\dots, \Gde(\zbf_N), \ \mubf) \\[0.5ex]
\text{s.t.} \quad
& \rbft_n(\zbf_0,\dots,\zbf_n, \ \mubf) = \mathbf{0}_{N_\zbf},
\qquad n = 0,1,\dots,N,\\
& \cbf(\Gde(\zbf_0),\dots, \Gde(\zbf_N), \ \mubf) \le \mathbf{0}_{N_\cbf}.
\end{aligned}
}
\label{eq:PDE_CO_LaSDI}
\end{equation}

We now turn to the computation of gradients required for gradient-based optimization.
The derivation closely parallels the full-order case presented in
Section~\ref{sec:gradient_optimzation}, with the key distinction that
sensitivities are computed with respect to the latent variables.
Defining $\ubft_n:= \Gde(\zbf_n) \in \mathbb{R}^{N_\ubf}$, the total derivative of the objective function with respect to parameter $\mu_i$ is
\begin{equation}
    \begin{aligned}
    \dd{f}{\mu_i}  (\ubft_0, \dots, \ubft_N, \ \mubf)
    =  & \sum_{n=0}^N \pd{f}{\ubf_n}  (\ubft_0, \dots, \ubft_N, \ \mubf)
    \ \nabla \Gde(\zbf_n)
    \  \pd{\zbf_n}{\mu_i}(\mubf) 
    +\pd{f}{\mu_i} (\ubft_0, \dots, \ubft_N, \ \mubf) . 
    \end{aligned}
    \label{eq:gradient_f_surrogate}
\end{equation}
Here and throughout, 
$\nabla \Gde(\zbf_n) \in \mathbb{R}^{N_\ubf \times N_\zbf}$ denotes the Jacobian of the decoder map $\Gde$ evaluated at $\zbf_n$.
The latent sensitivities $\pd{\zbf_n}{\mu_i}$ are obtained by differentiating $\rbft_n (\zbf_0,\dots,\zbf_n, \ \mubf),$
\begin{equation}
    \begin{aligned}
    %\dd{\rbf_0}{\mu_i}  &= \pd{\rbf_0}{\mu_i} + \pd{\rbf_0}{\ubf_0}\pd{\ubf_0}{\mu_i} = \mathbf{0}, \\
    \dd{\rbft_n}{\mu_i}  &= \pd{\rbft_n}{\mu_i}
    +
    \sum_{j=0}^n \pd{\rbft_n}{\zbf_{j}}  \ \pd{\zbf_j}{\mu_i}   = \mathbf{0}_{N_\zbf}, \qquad    n = 0, 1, \dots, N.
    \end{aligned}
    \label{eq:derivative_rtilden}
\end{equation}
The direct approach follows Eq.~\eqref{eq:direct}
\begin{equation}
    \begin{aligned}
        \pd{\zbf_0}{\mu_i} &= - \brackets{\pd{\rbft_0}{\zbf_0}}^{-1}\pd{\rbft_0}{\mu_i}, \\
        \pd{\zbf_{n}}{\mu_i} &=  - \brackets{\pd{\rbft_n}{\zbf_{n}}}^{-1}
        \parens{ \pd{\rbft_n}{\mu_i} +
        \sum_{j=0}^{n-1} \pd{\rbft_n}{\zbf_{j}} \ \pd{\zbf_{j}}{\mu_i}}, 
        \qquad n = 1, 2, \dots, N. 
    \end{aligned}
    \label{eq:direct_reduced}
\end{equation}
For the adjoint approach, we introduce adjoint variables $\ldbf_n \in \mathbb{R}^{N_\zbf}$ and follow the same procedure outlined in Eq.~\eqref{eq:adjoint_derivation}
\begin{equation*}
\begin{aligned}
    \sum_{n=0}^N \pd{f}{\ubf_n} 
    \ \nabla \Gde(\zbf_n) 
    \ \pd{\zbf_n}{\mu_i} = 
    -\sum_{n=0}^N \ldbf_n^T \ \pd{\rbft_n}{\mu_i}
    + \sum_{n=0}^N \parens{
    \pd{f}{\ubf_n}  \ \nabla \Gde(\zbf_n)  - \sum_{j=n}^N \ldbf_j^T \ \pd{\rbft_j}{\zbf_n}
    }\pd{\zbf_n}{\mu_i}.
\end{aligned}
\end{equation*}
Choosing the adjoint variables to eliminate the dependence on
$\pd{\zbf_n}{\mu_i}$ yields the backward-in-time adjoint equations
\begin{equation}
    \begin{aligned}
    \ldbf_N^T & = \pd{f}{\ubf_N} \  \nabla \Gde(\zbf_N) \ \brackets{\pd{\rbft_N}{\zbf_N}}^{-1},\\
        \ldbf_n^T  &= \parens{\pd{f}{\ubf_n} \ \nabla \Gde(\zbf_n) - \sum_{j=n+1}^N \ldbf_j^T \ \pd{\rbft_j}{\zbf_n}} 
        \brackets{\pd{\rbft_n}{\zbf_n}}^{-1}, 
        \qquad n = (N-1), (N-2), \dots, 0. 
    \end{aligned}
    \label{eq:adjoint_reduced}
\end{equation}

\begin{remark}
The computational advantages of LaSDI become apparent when comparing the cost of gradient evaluation.
In both the direct and adjoint formulations, the full-order model requires $N+1$ linear solves involving Jacobians
$\pd{\rbf_n}{\ubf_n} \in \Rbb^{N_\ubf \times N_\ubf}$.
In contrast, the reduced-order formulation involves $N+1$ linear solves with Jacobians
$\pd{\rbft_n}{\zbf_n}  \in \Rbb^{N_\zbf \times N_\zbf}$.
Since $N_\zbf \ll N_\ubf$, this reduction yields substantial savings in both memory usage and computational complexity.
\end{remark}

\begin{comment}
We begin by recalling the full-order model (FOM) describing the evolution of the state $\ubf(t, \ \mubf)$, governed by the semi-discrete system in Eq.~\eqref{eq:spatial_discretized_PDE}.
Within the LaSDI framework (see Section~\ref{sec:background}), the high-dimensional FOM state $\ubf \in \Rbb^{N_\ubf}$ is approximated by a low-dimensional latent state $\zbf \in \Rbb^{N_\zbf}$, with $N_\zbf << N_\ubf$ through encoder--decoder maps $\Gen$ and $\Gde$. The evolution of the latent state $\zbf(t, \ \mubf)$ is governed by the surrogate dynamics defined in Eq.~\eqref{eq:paramtric_latent_ode}. 
Applying a time-integration scheme to this latent dynamical system yields a corresponding set of \emph{reduced} residuals, denoted by $\rbft_n$, which parallel the full-order residuals $\rbf_n$ in Eq.~\eqref{eq:pde_residuals}. 
As in the full-order case, these residuals enforce the discrete-time evolution of the latent variables. For notational convenience, we write $\zbf_n \equiv \zbf_n(\mubf)$.
The reduced residuals are given by
\end{comment}

\section{Results}
\label{sec:results}
In this section, we evaluate the performance of WLaSDI surrogates for accelerating PDE-constrained optimization, with a focus on both accuracy and computational efficiency. We compare WLaSDI against optimization using full-order models (FOM) as well as alternative surrogate-based approaches, with particular attention to robustness under noisy training data. Specifically, we consider additive Gaussian noise applied to the snapshot data, with standard deviation $\sigma = \sigma_{NR}\norm{\Ubf}F$, where $\sigma_{NR}$ denotes the noise ratio.

We demonstrate the performance of WLaSDI on three benchmark problems: the one-dimensional inviscid Burgers' equation, the two-stream instability Vlasov-Poisson system, and thermal radiative transfer (TRT) for optimal hohlraum design. Across these examples, we consider both gradient-based and derivative-free optimization strategies, all implemented using \texttt{scipy.optimize} \cite{VirtanenGommersOliphantEtAl2020NatMethods}, and evaluate performance in terms of solution accuracy and wall-clock time.
Unless otherwise noted, all computations are performed on a single 2.0 GHz Intel Xeon CPU core of the LC Livermore Computing Dane system at LLNL, and reported runtimes are averaged over 10 runs.

\subsection{Burgers' Equation}
\label{sec:Burgers}

\subsubsection{High-Fidelity Simulation and WLaSDI Surrogates}
\label{sec:Burgers_FOM_ROM}

\begin{figure}[htbp]
\centering
    \includegraphics[width=0.4\textwidth]{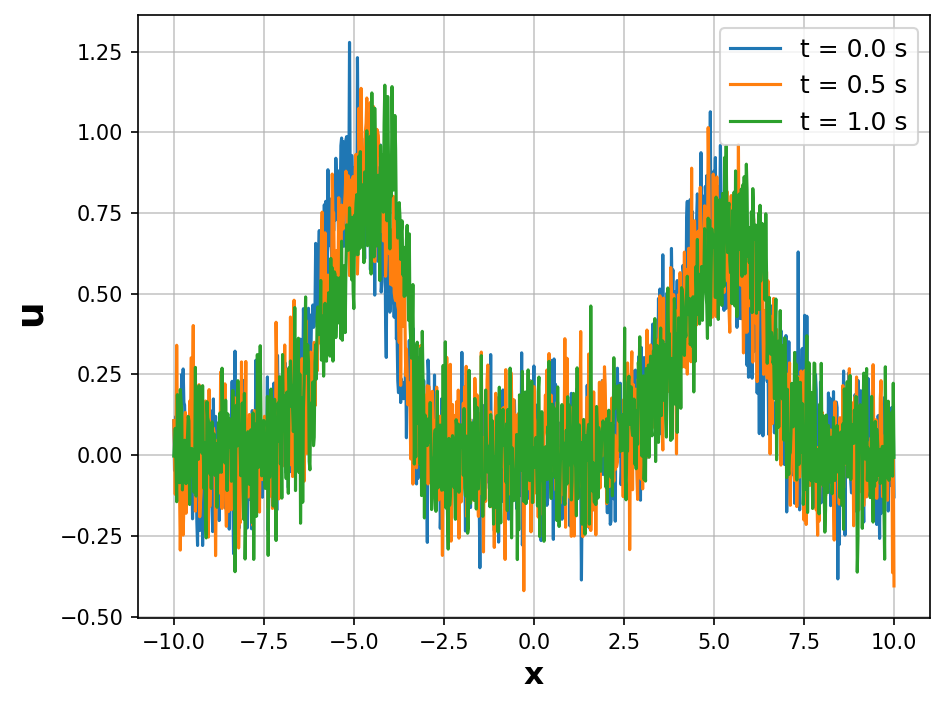}
%\caption{High-fidelity solution corresponding to $\mubf = [0.7, 1.1, 0.9, 0.9]$, with $10\%$ additive Gaussian noise.}
 \caption{\textbf{Burgers' Equation}: High-fidelity solution of Eq.~\eqref{eq:Burgers_eq} with $40\%$ noise used for training ($\mubf = [0.7, \, 1.1, \, 0.9, \, 0.9]$).}
\label{fig:Burgers_40_sim}
\end{figure}

We consider the one-dimensional inviscid Burgers' equation
\begin{equation}
\frac{\partial u}{\partial t} + u \ \frac{\partial u}{\partial x} = 0,  \qquad t \in [0, 1], \qquad  x \in [-10, 10],
\label{eq:Burgers_eq}
\end{equation}
subject to periodic boundary conditions  $u(x = 10, t) = u(x = -10, t)$. The initial condition is parameterized by $\mubf = [a_1, w_1, a_2, w_2] \in \Dcal$ and is defined as
\begin{equation}
u_0(x, \mubf) = u(x,\ t = 0; \ \mubf = [a_1, w_1, a_2, w_2]) = a_1\exp{\parens{-\frac{(x-5)^2}{2w_1^2}}} + a_2\exp{\parens{-\frac{(x+5)^2}{2w_2^2}}},
\label{eq:Burgers_IC}
\end{equation}
representing two Gaussian pulses with amplitudes $a_i$ and widths $w_i$. The parameter domain is 
$$\Dcal = [0.7, \, 0.9] \times [0.9, \, 1.1] \times [0.7, \, 0.9] \times [0.9, \, 1.1].$$
For numerical discretization, the spatial domain is approximated using a uniform grid with mesh spacing $\Delta x = 0.02$, and spatial derivatives are computed using a first-order forward finite difference scheme. Time integration is performed with the backward Euler method using a time step $\dt = 0.001$. Figure~\ref{fig:Burgers_40_sim} shows the numerical solution of Eq.~\eqref{eq:Burgers_eq} for $\mubf = [0.7, \, 1.1, \, 0.9, \, 0.9]$, with $40\%$ Gaussian noise.

To construct the WLaSDI surrogates, we follow the training procedure described in Section~\ref{sec:background}. The training dataset consists of 16 high-fidelity solutions of Eq.~\eqref{eq:Burgers_eq}, obtained from the parameterized initial condition in Eq.~\eqref{eq:Burgers_IC}. The parameters are sampled on a structured grid in the four-dimensional parameter space, using all combinations of
$$ a_1, a_2 \in \curlies{0.7, \, 0.9},  \quad w_1, w_2 \in \curlies{0.9, \, 1.1}.
$$
For compression, we employ Proper Orthogonal Decomposition (POD) with latent dimension $N_\zbf = 15$. The dimension is selected such that the first 15 singular values capture $99.99\%$ of the total singular value energy (defined as the cumulative sum of squared singular values divided by the total sum). To model the dynamics in the latent space, we adopt first-order polynomial features,
\begin{equation}
 \dd{\zbf}{t} (t) = \Wbf^T \thetabf \parens{\zbf(t) }, \quad \text{where} \quad  \thetabf(\zbf)  :=  \left[\begin{array}{ccccc}
1 & z_1 & z_2 \cdots &  z_{N_\zbf} \end{array}\right]^T \in \Rbb^{N_\zbf + 1}. 
\label{eq:latent_1st_order_poly}
\end{equation}
%as defined in Eq.~\eqref{eq:latent_1st_order_poly}.
%we record an order of magnutude speed up for rom
For the WLaSDI surrogate, the latent dynamics are identified using WENDy, as described in Section~\ref{sec:background}, together with the optimal test function construction introduced in \cite{TranBortz2026SISC}. We also employ the implicit parameterization strategy detailed in \ref{sec:implicit_parameterization}.
In addition to the noise-free setting, we consider noisy training data by introducing additive noise levels of $20\%$ and $40\%$.

\subsubsection{WLaSDI Performance in PDE-Constrained Optimization}
\label{sec:Burgers_optimization}

In this section, we investigate the performance of the WLaSDI surrogate in a PDE-constrained optimization setting. The objective is twofold: (i) to assess the robustness and accuracy of WLaSDI relative to other surrogate-based approaches, and (ii) to evaluate its computational efficiency compared to full-order model (FOM) optimization.

To this end, we compare optimization using WLaSDI-based surrogates against both alternative surrogate methods and direct optimization using the FOM. In the noise-free setting, we specifically benchmark WLaSDI against FOM-based optimization and radial basis function (RBF) interpolation of the objective function. To further assess robustness, we also consider surrogate models trained on noisy data, with noise levels of $20\%$ and $40\%$. The methods compared in this setting include WLaSDI, LaSDI \cite{FriesHeChoi2022ComputerMethodsinAppliedMechanicsandEngineering}, and RBF interpolation.

Optimization is performed using both derivative-free and gradient-based methods. For derivative-free optimization, we employ Constrained Optimization BY Quadratic Approximations (COBYQA) \cite{Ragonneau2022,RagonneauZhang2025}. For gradient-based optimization, we use the quasi-Newton Broyden-Fletcher-Goldfarb-Shanno (BFGS) algorithm \cite{NocedalWright2006}. Gradients are computed using adjoint-based techniques: the FOM gradients are obtained using Section~\ref{sec:gradient_optimzation}, while gradients for WLaSDI and LaSDI follow Section~\ref{sec:WLaSDI_optimization}.

\begin{figure}[htbp]
\centering
\includegraphics[width=0.4\textwidth]{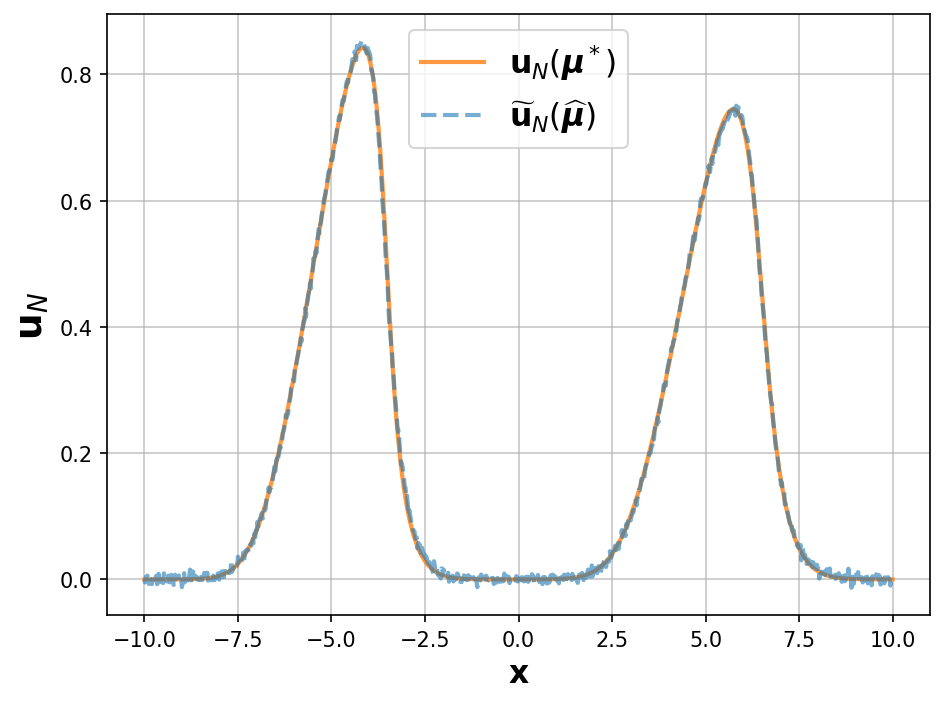}
\caption{\textbf{Burgers' Equation}: Comparison between the target state $\ubf_N(\mustar)$, (orange) and the recovered solution $\ubft_N(\muhat)$ (blue) obtained using WLaSDI model trained with $40\%$ noise. $\ubft_N(\muhat)$ closely matches the target state $\ubf_N(\mustar)$.}
\label{fig:Burgers_40_opt}  
\end{figure}

We consider an inverse problem in which the goal is to identify the parameter vector $\mubf$ defining the initial condition such that the solution at the final time $T = 1$ matches a prescribed target state. Specifically, 
\begin{equation}
    %f(\ubf_N, \mubf) = \norm{\ubf_N(\mubf) - \ubf_N(\mubf_{\mathrm{target}}) }_2
    f(\ubf_N (\mubf)) = \norm{\ubf_N(\mubf) - \ubf_N(\mubf^\star) }_2^2
    \label{eq:Burgers_obj}
\end{equation}
where $\ubf_N(\mubf)$ denotes the fully discrete solution at the final time corresponding to $\mubf$, and $\mustar$ is the target parameter vector. The target parameter is chosen as
$$\mustar = [0.75, \, 1.05, \, 0.85, \, 0.95],$$
with the associated target solution $\ubf_N(\mustar)$ shown in Fig.~\ref{fig:Burgers_40_opt} (orange).

Since the true solution of the inverse problem is known, we can directly quantify the accuracy of the recovered parameter $\muhat$. We evaluate performance using two metrics: (i) the relative parameter error,
\begin{equation}
    E_2(\muhat) := \frac{\norm{\muhat - \mustar}_2}{\norm{\mustar}_2},
    \label{eq:error_metric}
\end{equation}
and (ii) the true objective function evaluated using the FOM, which measures how well the recovered parameter reproduces the target state.
\begin{equation}
f(\ubf_N(\muhat)) := \norm{\ubf_N(\muhat) - \ubf_N(\mubf^\star)}_2^2.
\label{eq:obj_metric}
\end{equation}

Tables~\ref{tab:Burgers_0_noise}, \ref{tab:Burgers_20_noise}, and \ref{tab:Burgers_40_noise} present the optimization results for Burgers' equation using both derivative-free and gradient-based methods.
In Table~\ref{tab:Burgers_0_noise}, we compare WLaSDI-based optimization against full-order model (FOM) optimization and radial basis function (RBF) interpolation of the objective function. Tables~\ref{tab:Burgers_20_noise} and \ref{tab:Burgers_40_noise} focus on evaluating the robustness of surrogate-based methods when trained on noisy data.
In each table, we report the relative parameter error and the true objective function value, defined in Eqs.~\eqref{eq:error_metric} and \eqref{eq:obj_metric}, respectively. We also include the surrogate training time, as well as the optimization runtime for each method. In addition, the optimization cost is shown in terms of both runtime and the number of function evaluations.

\begin{table}[htbp]
\centering

\renewcommand{\arraystretch}{1.1}
\begin{tabular}{lccc}
\toprule
& WLaSDI & FOM & Interpolation \\
\midrule
$E_2(\muhat)$ (\%) 
& \num{ 0.4226923466}
& 0
& \num{6.7107824391} \\

$f(\ubf_N (\muhat))$  
& \num[round-mode=places,round-precision=4]{0.0015081896090065729}
&  0
&  \num{0.637454451051113} \\

Train time (s) 
& \num{25.579563856124878}
&  - & \num{0.028652191162109375} \\

\midrule
\multicolumn{4}{c}{\textbf{Derivative-free} (COBYQA)}\\
\midrule

Func.\ evals & 63 & 63 & 69 \\
Runtime (s) 
& \num{6.642259120941162}
& \num{112.0834424495697}
& \num{0.06338540712992351}\\

\midrule
\multicolumn{4}{c}{\textbf{Gradient-based} (BFGS)}\\
\midrule

Func.\ evals & 16 & 15 & -- \\
Grad.\ evals & 16 & 15 & -- \\
Runtime (s) 
& \num{2.5236499309539795} 
& \num{41.973286867141724}
&  - \\
\bottomrule
\end{tabular}

\caption{\textbf{Burgers' Equation}: Optimization results using noise-free surrogates and the full-order model (FOM), evaluated with both derivative-free (COBYQA) and gradient-based (BFGS) methods. WLaSDI outperforms RBF interpolation and achieves $20\times$ speedup over the FOM with high accuracy.}

\label{tab:Burgers_0_noise}
\end{table}

\begin{table}[htbp]
\centering

\renewcommand{\arraystretch}{1.1}
\begin{tabular}{lccc}
\toprule
& WLaSDI & LaSDI & Interpolation \\
\midrule

$E_2(\muhat)$ (\%) 
& \num{0.4731310328}
& \num{3.4810251782}
& \num{4.3778386364} \\

$f(\ubf_N (\muhat))$  
& \num[round-mode=places,round-precision=4]{0.0023319014400640012}
& \num{0.13995557364864564}
& \num{0.4370268725285183} \\

Train time (s) 
& \num{25.22743582725525} 
& \num{24.143343210220337}
& \num{0.03255939483642578} \\

\midrule
\multicolumn{4}{c}{\textbf{Derivative-free} (COBYQA)}\\
\midrule

Func.\ evals & 56 & 58 & 64\\

Runtime (s) 
& \num{4.99327921867370}
& \num{33.11495065689087}
& \num{0.060550530751546226}\\

\midrule
\multicolumn{4}{c}{\textbf{Gradient-based} (BFGS)}\\
\midrule

Func.\ evals & 14 & 16 & -- \\
Grad.\ evals & 14 & 16  & -- \\

Runtime (s) 
& \num{ 1.9426658153533936}
& \num{10.434960126876831} 
& - \\

\bottomrule
\end{tabular}

\caption{\textbf{Burgers' Equation}: Optimization results using surrogates trained with $20\%$ noise, evaluated with both derivative-free (COBYQA) and gradient-based (BFGS) methods. WLaSDI outperforms both LaSDI and RBF interpolation.}
\label{tab:Burgers_20_noise}
\end{table}

\begin{table}[htbp]
\centering

\renewcommand{\arraystretch}{1.1}
\begin{tabular}{lccc}
\toprule
& WLaSDI & LaSDI & Interpolation \\
\midrule
$E_2(\muhat)$ (\%) 
& \num{0.5065763165}
& \num{5.6929838584}
& \num{7.3716951683} \\

$f(\ubf_N (\muhat))$  
& \num[round-mode=places,round-precision=4]{0.0034995496033569913}
& \num{0.309040755185319}
& \num{1.4727747077161946} \\

Train time (s) 
& \num{25.22743582725525} 
& \num{24.143343210220337}
& \num{0.03255939483642578} \\

\midrule
\multicolumn{4}{c}{\textbf{Derivative-free} (COBYQA)}\\
\midrule

Func.\ evals & 56 & 58 & 57\\
Runtime (s) 
& \num{3.50203275680542}
& \num{29.631888389587402}
& \num{0.05662266413370768}\\

\midrule
\multicolumn{4}{c}{\textbf{Gradient-based} (BFGS)}\\
\midrule

Func.\ evals & 14 & 14 & --  \\
Grad.\ evals & 14 & 14 & --  \\
Runtime (s) 
& \num{1.4427263736724854}
& \num{9.651936531066895}
& - \\
\bottomrule
\end{tabular}
\caption{\textbf{Burgers' Equation}: Optimization results using surrogates trained with $40\%$ noise, evaluated with both derivative-free (COBYQA) and gradient-based (BFGS) methods. WLaSDI outperforms both LaSDI and RBF interpolation.}
\label{tab:Burgers_40_noise}
\end{table}

Across all scenarios, WLaSDI consistently outperforms the other surrogate approaches (LaSDI and RBF interpolation), achieving improvements of up to two orders of magnitude in accuracy, even under high noise levels. This robustness to noise is illustrated in Fig.~\ref{fig:Burgers_40_opt}, which compares the target state $\ubf_N(\mubf^\star)$ with the reconstructed solution $\ubft_N(\muhat)$ obtained using WLaSDI trained on $40\%$ noisy data (see also Fig.~\ref{fig:Burgers_40_sim}).

In terms of computational efficiency, WLaSDI provides a significant advantage over the FOM, achieving approximately a $20\times$ speedup, and about a $5\times$ speedup compared to LaSDI. Although WLaSDI is slower than interpolation-based methods, it yields substantially more accurate parameter estimates. Additionally, gradient-based optimization consistently converges faster than derivative-free optimization, particularly with adjoint-based implementation.

\subsection{Vlasov-Poisson Equation}
\label{sec:Vlasov}

\subsubsection{High-Fidelity Simulation and WLaSDI Surrogates}
\label{sec:Vlasov_FOM_ROM}

\begin{figure}[htbp]
    \centering
    \includegraphics[width=0.9\textwidth]{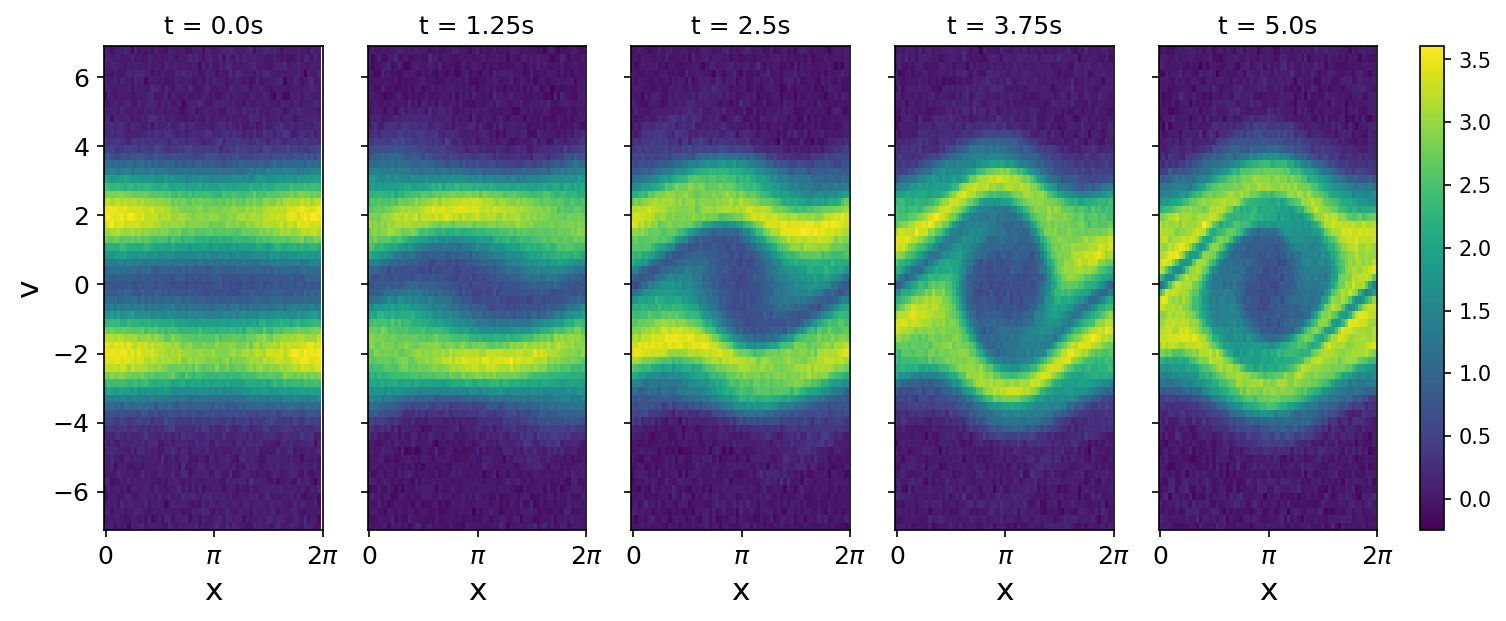}
    \caption{\textbf{Vlasov–Poisson:} 
    Phase-space evolution of the distribution function $u(x,v,t;\ \mubf)$ for $\mubf=[k, T]=[1, \, 0.9]$, including $5\%$ additive Gaussian noise.}
    \label{fig:Vlasov_simumation}
\end{figure}

We consider a simplified 1D-1V Vlasov-Poisson system describing the evolution of a collisionless plasma. The governing equations consist of a transport equation for the particle distribution function $u(x, v, t; \ \mubf)$ in phase space with periodic boundary conditions, coupled with Poisson's equation for the self-consistent electrostatic potential $\phi(x)$. The system is given by
\begin{equation}
    \begin{aligned}
        & \pd{u}{t} + v \pd{u}{x} + \pd{}{v} \! \brackets{\dd{\phi}{x} \ u} = 0, \quad t \in [0, 5], \quad x \in [0, 2\pi], \quad v \in [-7, 7],\\
        & \frac{d^2 \phi}{dx^2} = \int_v u dv.
    \end{aligned}
    \label{eq:Vlasov_Poisson}
\end{equation}
We focus on the two-stream instability case,  where the initial condition consists of two counter-streaming Maxwellian beams centered at velocities $v \pm 2 $. 
This initial condition is parameterized by the temperature $T$ and the spatial perturbation wave number $k$, $\mubf = [k, T]$. The initial condition is given by 
\begin{equation*}
    u(x, v, t = 0; \ \mubf = [k, T]) = \frac{8}{\sqrt{2 \pi T}} \parens{
    1 + \frac{1}{10} \cos(kx)}
    \parens{\exp \parens{ -\frac{(v-2)^2}{2T}} +  \exp \parens{-\frac{(v+2)^2}{2T}}}.
\end{equation*}
An example of this initial condition for $\mubf = [1, \, 0.9]$ is shown in the leftmost panel of Fig.~\ref{fig:Vlasov_simumation}. 
The parameter domain considered here is $\Dcal =  [1.0, \, 1.2] \times [0.9, \, 1.1]$.

High-fidelity solutions are computed using the HyPar solver \cite{GhoshLoffeldKim}. Spatial discretization is performed with a weighted essentially non-oscillatory (WENO) finite difference scheme \cite{JiangShu1996JournalofComputationalPhysics},
and time integration is carried out using a fourth-order explicit Runge–Kutta method. A representative simulation corresponding to $\mubf = [1, 0.9]$, with $5\%$ additive Gaussian white noise, is shown in Fig.~\ref{fig:Vlasov_simumation}. 

\begin{figure}[htbp]
    \centering
    \includegraphics[width=0.4\textwidth]{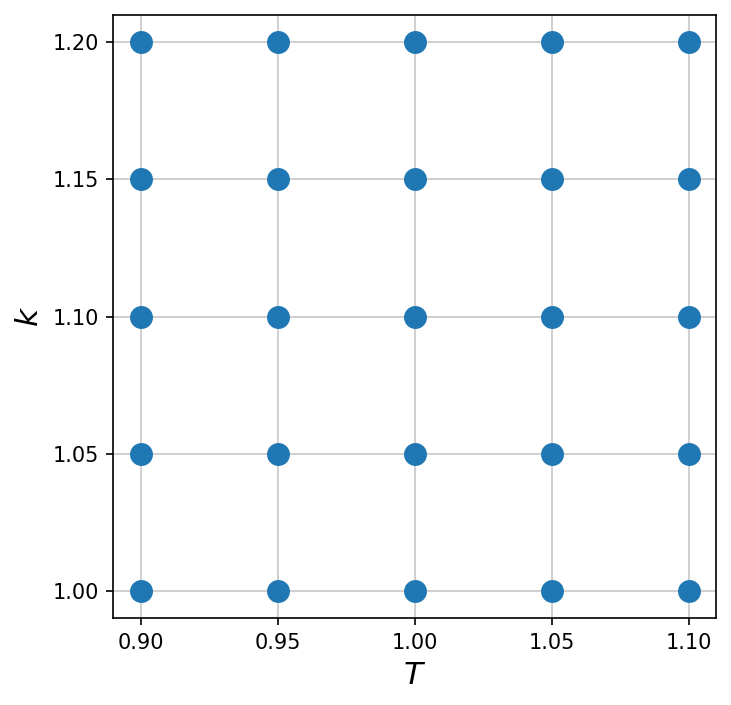}
    \caption{\textbf{Vlasov-Poisson:} Parameter space  $\Dcal =  [1.0, \, 1.2] \times [0.9, \, 1.1]$, with training parameters shown in blue.}
    \label{fig:Vlasov_training_params}
\end{figure}

To construct the WLaSDI surrogates, we employ $25$ training parameter samples arranged on a uniform grid in parameter space, as shown in Fig.~\ref{fig:Vlasov_training_params}. Two surrogate models are constructed: one trained on noise-free data and another trained on data with $5\%$ noise.
These models are trained using the autoencoder framework introduced in \cite{HeTranBortzEtAl2025NumericalMethEngineering}, where the autoencoder and the latent ODE are learned simultaneously. Each autoencoder consists of three hidden layers with $1000$, $100$, and $50$ neurons, respectively, using ReLU activation functions. A latent dimension of $6$ is used in both cases.
The latent dynamics are modeled using a first-order polynomial representation, as described in Eq.~\eqref{eq:latent_1st_order_poly}. 
%The coefficient matrices $\curlies{\Wbf^{(k)}}$ are determined using WENDy (Section~\ref{sec:WENDy}), with optimal test function hyperparameters selected according to \cite{TranBortz2026SISC}. 
Gaussian process regression (\ref{sec:interp_gp}) is used to interpolate the coefficient matrices across the parameter space.
For a detailed evaluation of the WLaSDI method in terms of simulation error for this problem, we refer the reader to \cite{HeTranBortzEtAl2025NumericalMethEngineering}.

In terms of computational cost, a single full-order model (FOM) simulation requires approximately $72$ seconds on $16$ cores, whereas a single WLaSDI evaluation takes $\num[round-mode=places,round-precision=1]{0.08712606430053711}$ seconds on a single core. This corresponds to a speedup of approximately five orders of magnitude.

\subsubsection{WLaSDI Performance in PDE-Constrained Optimization}
\label{sec:Vlasov_optimization}

\begin{figure}[htbp]
    \centering
    \includegraphics[width=0.75\textwidth]
    {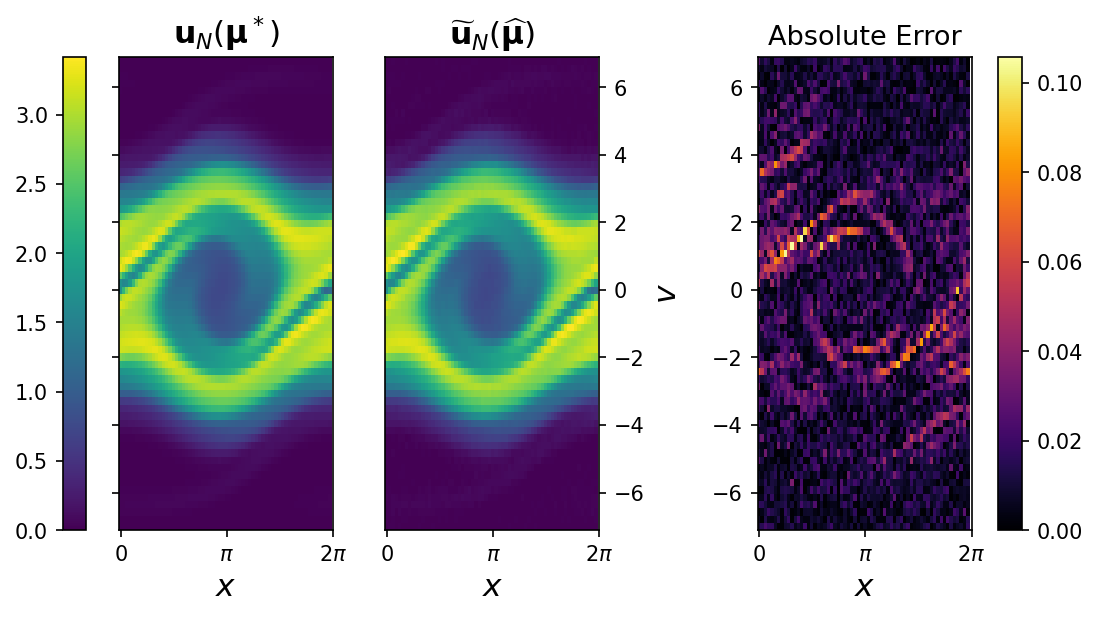}
 %\caption{\textbf{Vlasov-Poisson}: Target solution $\ubf_N(\mubf^\star)$, for $\mubf^\star = [1.06, 0.97]$; reconstructed solution $\ubft_N(\widehat{\mubf})$ from WLaSDI surrogate with $5\%$ noise, $\widehat{\mubf} = [\num[round-mode=places,round-precision=3]{1.06371124},\num[round-mode=places,round-precision=3]{0.95917576}]$; and the corresponding absolute error.}
    \caption{\textbf{Vlasov-Poisson}: 
        Target solution $\ubf_N(\mubf^\star = [1.16, \,  1.03])$;
        reconstructed solution $\ubft_N(\muhat =  [\num[round-mode=places,round-precision=3]{1.15875343}, \,
\num[round-mode=places,round-precision=3]{1.02288916}])$ from WLaSDI surrogate with $5\%$ noise,
and the absolute error $|\ubf_N(\mustar) -  \ubft_N(\muhat)|$.}
     \label{fig:Vlasov_opt_5}
\end{figure}

To demonstrate the performance of the WLaSDI surrogate in a PDE-constrained optimization setting, we consider an inverse problem in which the goal is to identify the parameter vector $\mubf$ defining the initial condition such that the solution at the final time $T = 5$ matches a prescribed target state.
Let $\ubf_N(\mubf)$ denote the fully discrete solution to the Vlasov-Poisson system (Eq.~\eqref{eq:Vlasov_Poisson}) with parameters $\mubf$ at the final time $t_N$. The objective function is defined as
\begin{equation}
    %f(\ubf_N, \mubf) = \norm{\ubf_N(\mubf) - \ubf_N(\mubf_{\mathrm{target}}) }_2
    f(\ubf_N, \mubf) = \norm{\ubf_N(\mubf) - \ubf_N(\mubf^\star) }_2^2
    \label{eq:Vlasov_obj}
\end{equation}
where $\mubf^\star$ denotes the parameter vector associated with the desired solution.
Here, the target parameter is chosen as $\mubf^\star = [1.16, \, 1.03]$. The corresponding target solution at $T = 5$ is shown in the leftmost panels of Fig.~\ref{fig:Vlasov_opt_5}.

Tables~\ref{tab:Vlasov_0_noise} and \ref{tab:Vlasov_5_noise} summarize the results of the optimization problem defined in Eq.~\eqref{eq:Vlasov_obj} for the noise-free and $5\%$ noisy training cases, respectively. In each table, we report results obtained using both a derivative-free method (COBYQA) and a gradient-based method (BFGS). To assess performance, we also report the error and the true objective value at the solution, as defined in Eqs.~\eqref{eq:error_metric} and \eqref{eq:obj_metric}.
In addition, we provide a rough estimate of the computational cost associated with the full-order model (FOM). The estimated FOM runtime is based on the total number of objective evaluations required during optimization, where each gradient evaluation, approximated via finite differences, is counted as two objective evaluations.

\begin{table}[htbp]
\centering

\begin{subtable}{\textwidth}
\centering
\renewcommand{\arraystretch}{1.1}
\begin{tabular}{lcc}
\toprule
& \textbf{Derivative-free} & \textbf{Gradient-based} \\
& (COBYQA) & (BFGS) \\
\midrule

$E_2(\muhat)$ (\%) 
& \num{0.23019240}
& \num{0.2301869381}\\

$f(\ubf_N (\muhat))$  
& $\num[round-mode=places,round-precision=3]{0.24056582814651217}$
& $\num[round-mode=places,round-precision=3]{0.24056582814651217}$\\

\midrule

Func. evals 
& $29$
& $11$ \\

Grad. evals 
& -- 
& $11$ \\

WLaSDI runtime (s)
& $\num[round-mode=places,round-precision=2]{0.64874267578125}$
& $\num[round-mode=places,round-precision=2]{3.6721086502075195}$ \\

\midrule 
Est. FOM cost (s)
& $\num[round-mode=places,round-precision=0]{33408}$
& $\num[round-mode=places,round-precision=0]{25344}$ \\

\bottomrule
\end{tabular}
\caption{Noise-free.}
\label{tab:Vlasov_0_noise}
\end{subtable}

\vspace{0.5cm}

\begin{subtable}{\textwidth}
\centering
\renewcommand{\arraystretch}{1.1}
\begin{tabular}{lcc}
\toprule
& \textbf{Derivative-free} & \textbf{Gradient-based} \\
& (COBYQA) & (BFGS) \\
\midrule

$E_2(\muhat)$ (\%) 
& \num{0.44961133}
& \num{ 0.4496257889}\\

$f(\ubf_N (\muhat))$  
& $\num[round-mode=places,round-precision=3]{2.365434673725959}$
& $\num[round-mode=places,round-precision=3]{2.365434673725959}$\\

\midrule

Func. evals 
& $32$
& $12$ \\

Grad. evals 
& -- 
& $12$ \\

WLaSDI runtime (s)
& $\num[round-mode=places,round-precision=2]{1.2841849327087402}$
& $\num[round-mode=places,round-precision=2]{4.392463684082031}$ \\
\bottomrule
\end{tabular}
\caption{$5\%$ noise.}
\label{tab:Vlasov_5_noise}
\end{subtable}

\caption{\textbf{Vlasov-Poisson}: Optimization results using WLaSDI surrogates. 
Results are presented for both the derivative-free (COBYQA) and gradient-based (BFGS) methods. (a) Surrogates trained with noise-free data. 
(b) Surrogates trained with $5\%$ noise.}
\label{tab:Vlasov}
\end{table}

In all cases, the recovered $\muhat$ are very close to the target $\mustar$. Figure~\ref{fig:Vlasov_opt_5} compares the target solution $\ubf_N(\mubf^\star)$ with the optimized solution $\ubft_N(\muhat)$ obtained using the WLaSDI surrogate trained on $5\%$ noisy data, together with the corresponding absolute error $|\ubf_N(\mustar) -  \ubft_N(\muhat)|$. The optimized solution closely matches the target state, demonstrating the robustness of WLaSDI. 

As reported in Tables~\ref{tab:Vlasov_0_noise} and \ref{tab:Vlasov_5_noise}, the gradient-based method requires fewer function (and gradient) evaluations than the derivative-free approach. However, the total WLaSDI runtime of the derivative-free method is lower than that of the gradient-based method.
This behavior can be explained by the higher cost of computing gradients. In particular, evaluating the gradient $\dd{f}{\mu_i}$ is significantly more expensive than evaluating the objective function $f$, with each gradient evaluation being approximately an order of magnitude more costly. This additional expense arises from differentiating through the nonlinear autoencoder and the interpolation of the latent space coefficient matrices, i.e., computing $\pd{\Wbf}{\mu_i}$ (see \ref{sec:RK}).
Moreover, since the parameter space is only two-dimensional, the number of function evaluations required by the derivative-free method remains relatively small. As a result, the computational advantage of the gradient-based approach is not observed in this setting. Nevertheless, WLaSDI-based optimization remains substantially faster than solving the problem using the full-order model.

\subsection{Thermal Radiative Transfer} 
\label{sec:TRT}
Nuclear fusion has long been pursued as a potential source of clean and renewable energy. In indirect-drive fusion experiments, a small fuel capsule is placed inside a device called a \emph{hohlraum}, where intense radiation is used to compress the capsule and initiate fusion reactions. For this process to succeed, the capsule must be heated as uniformly as possible; even small asymmetries can disrupt the compression and reduce performance. Achieving such uniform energy delivery is therefore a central and challenging design goal.

In this section, we formulate an optimization problem to identify temperature configurations that promote uniform heating of the fusion capsule. Solving this problem directly using the underlying high-fidelity model is computationally expensive. To address this challenge, we employ the WLaSDI surrogates to enable efficient exploration of the design space.

\subsubsection{Governing Equations and High-Fidelity Simulation}
\label{sec:TRT_FOM}
\begin{figure}[htbp]
    \centering
    \begin{subfigure}{0.49\textwidth}
        \centering
        \includegraphics[width=\textwidth]{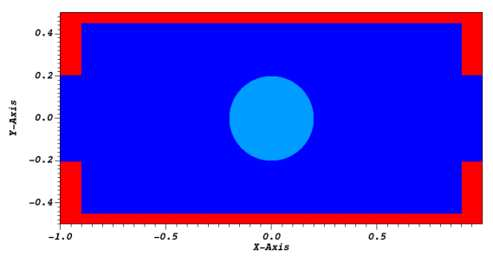}
         \caption{Idealized 2D hohlraum configuration showing the outer wall (red), gas-filled cavity (dark blue), and circular fuel capsule (light blue).}
         \label{fig:hohlraum_materials}
    \end{subfigure}
    \begin{subfigure}{0.49\textwidth}
        \centering
       \includegraphics[width=\textwidth]{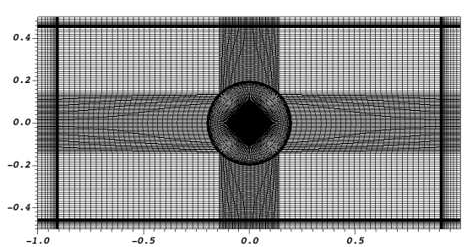}
        \caption{Spatial mesh consisting of $28,500$ elements used for the Discontinuous Galerkin (DG) discretization.}
        \label{fig:TRT_mesh}
    \end{subfigure} 
    \caption{\textbf{Thermal Radiative Transfer}: Computational setup for the hohlraum simulation: (a) Idealized hohlraum geometry and material regions; (b) DG spatial mesh.}
 \label{fig:TRT_setup}
\end{figure}

We consider the two-dimensional nonlinear Thermal Radiative Transfer (TRT) equations in an idealized 2D hohlraum configuration. The computational domain represents a simplified hohlraum chamber, illustrated in Fig.~\ref{fig:hohlraum_materials}. The outer wall is shown in red, the gas-filled cavity in dark blue, and the circular capsule in light blue. The spatial coordinate is denoted by $\xbf = (x,y)\in \Xcal$, where $\Xcal := [-1, \, 1] \times [-0.5, \, 0.5]$. 
The fuel capsule is centered at the origin $(0, 0)$, with the prescribed radius $r_0 = 0.2$.

To model thermal radiation transport within this configuration, we consider the frequency-averaged, grey TRT equations. These equations describe the coupled evolution of the specific intensity (angular flux) $\psi (\xbf, \omegabf, t)$ and the material temperature  $T (\xbf, t)$. The angular variable satisfies $\omegabf \in \mathbb{S}^2$ (the unit sphere), and time $t \in [0,t_{\mathrm{final}}]$. The governing equations are given by
\begin{equation}
    \begin{aligned}
 \frac{1}{c} \pd{\psi}{t} + \omegabf \cdot \nabla_{\xbf}\psi 
       & = \sigma(T)\parens{\frac{ac}{4\pi}T^{4}-\psi},\\
\psi(\xbf, \omegabf, t) & = \psi^{\mathrm{inc}}(\xbf, \omegabf, t), \quad  \omegabf \cdot \nbf(\xbf) \leq 0, \quad \xbf \in \partial \Dcal.
    \end{aligned}
    \label{eq:TRT_flux}
\end{equation}
\begin{equation}
    \begin{aligned}
 \rho c_{v}(T) \ \pd{T}{t} &= \sigma(T)  \int_{\mathbb{S}^{2}} \parens{ \psi - \frac{ac}{4\pi} T^4 } d \omegabf, \\
 T(\xbf, 0) & = T_0(\xbf).
    \end{aligned}
    \label{eq:TRT_temp}
\end{equation}
Here, $\nbf(\xbf)$ denotes the outward unit normal vector on the spatial boundary $\partial \Dcal$, and $\psi^{\mathrm{inc}}$ represents the prescribed inflow radiation. The material density is denoted by $\rho$, $c_v(T)$ is the temperature-dependent heat capacity, $c$ denotes the speed of light, and $a$ denotes the radiation constant.

We use microseconds, cm, EU, and keV as the units for time, length, energy, and temperature, respectively, where $\ensuremath{1\,\text{EU}=10^{12}\,\text{ergs}}$. In these units, the radiation constant and speed of light are given by $\ensuremath{a=137.199\,\text{cm}^{4}/\text{keV}^{4}}$ and $\ensuremath{c=2.9979\times10^{4}\,\text{cm}/\mu\text{s}}$. The absorption opacity $\sigma$ has units of $\text{cm}^{-1}$, and $\rho c_v$ has units of $\text{EU}/\text{cm}^3/\text{keV}$.

A fixed time step of $\ensuremath{0.0625\times10^{-5}}$ was used to resolve the relevant time scales.

As illustrated in Fig.~\ref{fig:hohlraum_materials}, the hohlraum geometry consists of three distinct material regions: an optically thick outer wall in red, an optically thin gas-filled cavity in dark blue, and a circular capsule in light blue. The material properties for each region are summarized in Table~\ref{tab:hohlraum_materials}.
\begin{table}[htbp]
\centering
\renewcommand{\arraystretch}{1.1}
\begin{tabular}{lcc}
\toprule
Region & Absorption opacity $\sigma$ & Heat capacity $\rho c_v$ \\
\midrule
Outer wall (red)      & $10^{4}$      & $10^{3}$ \\
Gas-filled cavity (dark blue) & $10^{-3}$     & $10^{-4}$ \\
Circular capsule (light blue) & $15$          & $2\times10^{-2}$ \\
\bottomrule
\end{tabular}
\caption{\textbf{Thermal Radiative Transfer}: Material properties for the three regions of the hohlraum configuration.}
\label{tab:hohlraum_materials}
\end{table}

\begin{figure}[htbp]
    \centering
    \includegraphics[width=0.9\textwidth]{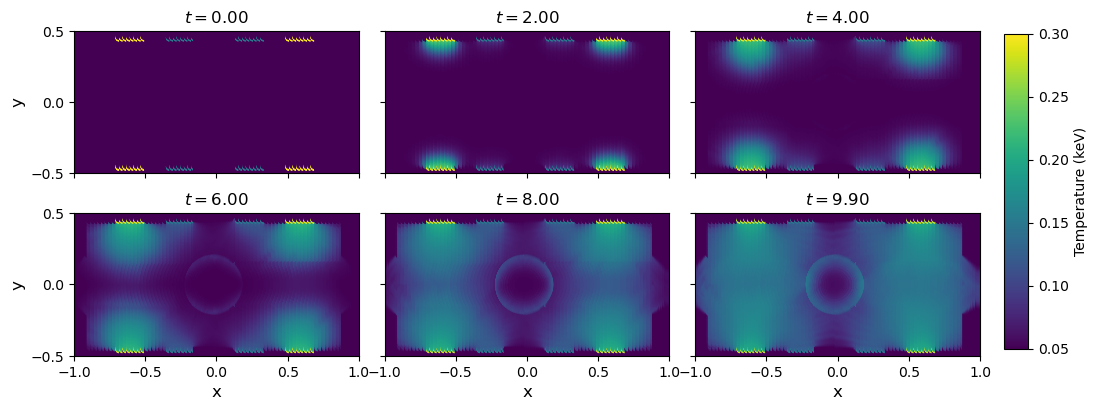}
    \caption{\textbf{Thermal Radiative Transfer}: Time evolution of the temperature field $T(\xbf,t;\mubf)$ for $\mubf =[T_{\mathrm{outer}},\,T_{\mathrm{inner}}] = [0.3,\,0.15]$ keV. 
    Corner and central hot spots are initialized at $T_{\mathrm{outer}}$ and  $T_{\mathrm{inner}}$, respectively, while the remainder of the domain is initialized at $0.05$ keV.}
    \label{fig:TRT_simulation}
\end{figure}

The temperature is initialized to a uniform background value of $0.05$ keV throughout the domain, except at eight
hot spots along the walls, as shown in the upper-left panel of Fig.~\ref{fig:TRT_simulation}. 
The initial condition is defined as a parameterized function $T_0(\xbf, \mubf)$, where
\begin{equation*}
    \mubf = [T_{\mathrm{outer}},\, T_{\mathrm{inner}}] \in [0.05,\,0.3]^2 \ \mathrm{(keV)}.
\end{equation*}
Specifically, $T_{\mathrm{outer}}$ denotes the temperature assigned to the four corner hot spots, while $T_{\mathrm{inner}}$ specifies the temperature of the four hot spots located closer to the center of the walls.
These two temperatures define the parameter space over which the reduced model and subsequent design optimization are performed.
Outside of these prescribed locations, the temperature remains at the uniform background value of $0.05$ keV. An example corresponding to $\mubf = [0.3, \, 0.15]$ keV is shown in Fig.~\ref{fig:TRT_simulation}.

The high-fidelity model is constructed by discretizing the angular domain using the discrete ordinates $S_N$ method \cite{CarlsonLathropLaboratory1965, LarsenMorel2010NuclearComputationalScience}. An angular product quadrature of order $30$ is employed, yielding $1,800$ discrete angular directions. The spatial domain is discretized using a Discontinuous Galerkin (DG) formulation with linear basis functions and standard upwinding on a mesh consisting of $28,500$ elements (see Fig.~\ref{fig:TRT_mesh}). Time integration is performed using a backward Euler scheme. The resulting nonlinear implicit systems at each time step are efficiently solved using the Variable Eddington Factor (VEF) method \cite{OlivierPaznerHautEtAl2023JournalofComputationalPhysics}. A single high-fidelity simulation requires approximately 30 minutes when executed on 72 CPU cores in parallel.

For this example, WLaSDI is used to simulate the parameterized temperature field $T(\xbf, t; \, \mubf)$ of the hohlraum device. For consistency with the design optimization detailed below, we use the notation $T(\xbf, t; \, \mubf)$ and $u(\xbf, t;  \, \mubf)$ interchangeably to denote the temperature field. Figure~\ref{fig:TRT_simulation} illustrates the time evolution of the temperature field corresponding to  $\mubf = [0.3, \, 0.15]$ keV.

\subsubsection{WLaSDI Surrogate Construction}
\label{sec:TRT_WLaSDI}

%and evaluate their predictive performance for the Thermal Radiative Transfer (TRT) system.
%In addition to the weak form formulation, we also construct surrogate models using the strong form (LaSDI) in order to assess the relative robustness of the two approaches to noisy training data. In particular, we consider both noise-free training data and data contaminated with $5\%$ additive Gaussian white noise.

\begin{figure}[htbp]
    \centering
    \includegraphics[width=0.4\textwidth]{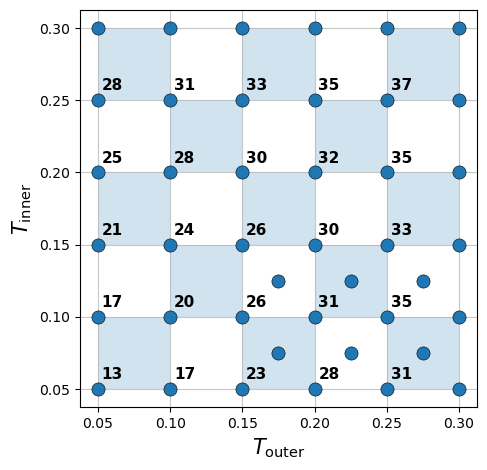}
\caption{\textbf{Thermal Radiative Transfer}:
    Decomposition of the training parameter space. 
    Each checkerboard cell defines a local WLaSDI submodel constructed from the enclosed training parameter (blue). The latent dimension $N_{\zbf}$ of each submodel is shown.}
    \label{fig:TRT_training_params}
\end{figure}

In this section, we describe the construction of WLaSDI surrogate models for both noise-free and noisy training data.
We first build the compression and decompression operators as described in Section~\ref{sec:background}. The two-dimensional parameter space $(T_{\mathrm{outer}},\, T_{\mathrm{inner}})$ is discretized using a uniform $5 \times 5$ grid with spacing $0.05$, as illustrated in Fig.~\ref{fig:TRT_training_params}. For each grid location, we construct a local surrogate submodel using a neighborhood consisting of $4$--$5$ enclosed parameter points. While this partitioning is not necessarily optimal in terms of training data sampling or parameter-space decomposition, and the development of more effective strategies is left for future work, it enables parallelized training of the local compression and decompression operators, reduces computational cost, and better captures variations in local solution complexity across the parameter domain.

For each local surrogate, Proper Orthogonal Decomposition (POD) is employed as the compression step. The solution snapshots are assembled into a data matrix, and singular value decomposition (SVD) is performed. The latent space dimension $N_\zbf$ is selected such that the retained singular values capture $99.9\%$ of the total singular value energy, defined as the ratio of the cumulative sum of the retained singular values to the total sum of all singular values. Across the $25$ local models, the resulting latent dimensions range from $13$ to $37$, as illustrated in Fig.~\ref{fig:TRT_training_params}.
To model the dynamics in the latent space, we adopt first-order polynomial features. The resulting latent ODEs take the form
\begin{equation*}
 \dd{\zbf}{t} (t) = \Wbf^T \thetabf \parens{\zbf(t) }, \quad \text{where} \quad  \thetabf(\zbf)  :=  \left[\begin{array}{ccccc}
1 & z_1 & z_2 \cdots &  z_{N_\zbf} \end{array}\right]^T \in \Rbb^{N_\zbf + 1}. 
%\label{eq:latent_1st_order_poly}
\end{equation*}
For each local WLaSDI submodel, the coefficient matrix $\Wbf$ is computed using the weak-form identification procedure described in Section~\ref{sec:WENDy} and further examined in \cite{TranHeMessengerEtAl2024ComputerMethodsinAppliedMechanicsandEngineering,HeTranBortzEtAl2025NumericalMethEngineering}. In addition, we incorporate the test function construction method introduced in \cite{TranBortz2026SISC} to further improve the accuracy of latent dynamics identification. 
For comparison, we also construct strong-form LaSDI models in which the coefficient matrix $\Wbf$ is identified using the SINDy framework \cite{BruntonProctorKutz2016ProcNatlAcadSciUSA}. For both approaches, we employ the global parameterization strategy described in \ref{sec:global}.

WLaSDI  provides substantial computational acceleration. The total training time for the local WLaSDI surrogate models is $\num[round-mode=places,round-precision=0]{2271.956058}$ seconds. Once trained, the simulation time per parameter instance is approximately $2.76$ seconds on a single CPU core.
Compared to a single high-fidelity TRT simulation, this corresponds to a wall-clock speed-up of approximately $6.5 \times 10^{2}$ and a one-core equivalent CPU-time speed-up of approximately $4.5 \times 10^{4}$ per evaluation. Thus, WLaSDI reduces the computational cost by more than two orders of magnitude in wall-clock time and over four orders of magnitude in CPU time. These substantial gains make the WLaSDI surrogate model well-suited for efficient and accurate design optimization of the hohlraum device detailed below.

\subsubsection{Design Optimization of the Hohlraum}
\label{sec:TRT_optimization}

To achieve thermonuclear fusion, the fuel capsule located at the center of the hohlraum (see Fig.~\ref{fig:TRT_setup}) must be heated as uniformly as possible. In this work, we seek to promote uniform heating by adjusting the prescribed temperatures of the laser-induced wall hot spots. 
Therefore, the quantity of interest is the total temporal variation in the temperature distribution along the surface of the capsule. Specifically, we measure the standard deviation of the capsule surface temperature at each time step and integrate this quantity over time. 

The temperature along the capsule surface is obtained by linearly interpolating the Discontinuous Galerkin (DG) elements onto points along a circular contour centered at the origin with the prescribed capsule radius. 
Let $\ubf_n(\mubf)$ denote the discrete temperature field at time step $t_n$, we obtain an interpolation operator
$ \mathcal{I}_\Gamma$, which maps the DG solution to temperature values along the capsule surface $\Gamma$, 
where $\Gamma$ denotes the circular boundary of the capsule.

Let $\sigma(\cdot)$ denote the standard deviation of the interpolated surface temperature. To promote uniform heating, we define the objective function as the time-integrated standard deviation,
\begin{equation}
    f(\ubf_0, \cdots, \ubf_N; \ \mubf) = \dt  \sum_{n=0}^{N}{}'' \sigma \parens{\mathcal{I}_\Gamma \ubf_n}
    \label{eq:TRT_objective}
\end{equation} 
where the double prime notation on the sum $\sum{}''$ denotes the trapezoidal-weighted sum, i.e., the first and last terms are halved before summing. 

\begin{figure}[htbp]
    \centering
    \includegraphics[width=0.4\linewidth]{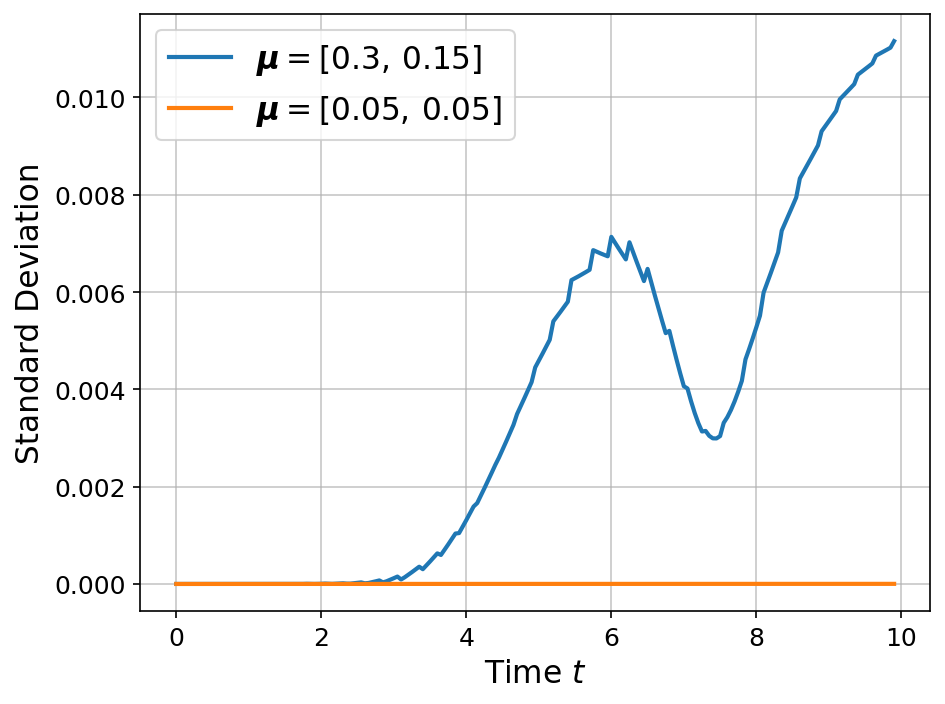}
    \caption{\textbf{Thermal Radiative Transfer:} Standard deviation of the capsule surface temperature for $\mubf = [0.3,\,0.15]$ (blue) and $\mubf = [0.05,\,0.05]$ (orange). The standard deviation for $\mubf = [0.05,\,0.05]$ is negligible, indicating nearly uniform temperature along the capsule surface.}
    \label{fig:TRT_QoI_0}
\end{figure}

Figure~\ref{fig:TRT_QoI_0} illustrates the time evolution of $\sigma \parens{\mathcal{I}_\Gamma \ubf_n}$ for two parameter choices: $\mubf = [0.3, \, 0.15]$ (blue) and $\mubf = [0.05,\, 0.05]$ (orange). The latter yields an almost negligible standard deviation, indicating nearly uniform temperature throughout the simulation. This behavior is expected, as the initial temperature is $0.05$ keV and the hot-spot temperatures are set to the same value, resulting in minimal thermal gradients and thus little evolution. While this represents a trivial optimum of the objective, it is not physically interesting.

To ensure a physically meaningful design, we impose an additional constraint requiring the mean capsule surface temperature at the final time to exceed a prescribed threshold $T_{\mathrm{min}} = 0.1$ keV. Let $\mu(\cdot)$ denote the mean operator,  the constraint is defined as
\begin{equation}
    g(\ubf_N; \ \mubf) = \mu \parens{\mathcal{I}_\Gamma \ubf_N} - T_{\mathrm{min}} \geq 0.
    \label{eq:TRT_constraint}
\end{equation} 
We now seek to minimize the objective functional in Eq.~\eqref{eq:TRT_objective}, subject to the governing TRT equations (Eqs.~\eqref{eq:TRT_flux}–\eqref{eq:TRT_temp}) and the inequality constraint in Eq.~\eqref{eq:TRT_constraint}. The optimization is carried out using pretrained WLaSDI surrogate models, as direct optimization with the full-order model is computationally prohibitive. To assess robustness to noise, we additionally compare against LaSDI surrogate models.

Because the WLaSDI/LaSDI surrogate models consist of 25 local submodels (Fig.~\ref{fig:TRT_training_params}), the resulting objective function is piecewise-defined and may be discontinuous along the boundaries between subregions. In this setting, global derivative-free optimization methods such as Differential Evolution \cite{StornPrice1997JournalofGlobalOptimization} can be applied directly. However, gradient-based methods are typically much more computationally efficient, although the piecewise structure of the surrogate means that the objective is discontinuous at subregion interfaces.
One natural strategy is to perform the optimization independently on each submodel and then combine the results. This approach is highly parallelizable and allows the optimization tasks to be distributed across multiple processors. Alternatively, one may perform optimization directly on the full surrogate model. In this case, the objective is piecewise defined and may be nondifferentiable at the boundaries between subregions. A common remedy would be to smooth the objective to obtain a continuous optimization problem. In this case, however, the discontinuities occur only along subregion boundaries, which form a set of measure zero in the parameter space. Consequently, the gradient is well defined \emph{almost everywhere}, allowing gradient-based optimization methods to be applied effectively in practice.

In the following, we present results for three optimization strategies: 
\begin{itemize}
\item  \hyperref[par:TRT_DiffEv]{Global derivative-free optimization} using Differential Evolution \cite{StornPrice1997JournalofGlobalOptimization}.
\item \hyperref[par:TRT_parallel]{Parallel gradient-based optimization across the $25$ local submodels} using  Sequential Least Squares Programming (SLSQP) \cite{LawsonHanson1995}.
\item \hyperref[par:TRT_direct]{Direct gradient-based optimization} of the full surrogate model using trust-region constrained \cite{ByrdHribarNocedal1999SIAMJOptim}.
\end{itemize} 
For each strategy, we compare results obtained using WLaSDI and LaSDI surrogates.

\paragraph{Global derivative-free optimization}
\label{par:TRT_DiffEv}

We first present results for global optimization using the derivative-free algorithm Differential Evolution \cite{StornPrice1997JournalofGlobalOptimization}. The resulting global optimum is then used as a baseline to evaluate the performance of gradient-based optimization methods.

Table~\ref{tab:TRT_DiffEv} summarizes the optimization results obtained using surrogate models trained with noise-free data and with $5\%$ noisy training data. In the noise-free setting, the LaSDI and WLaSDI surrogates produce nearly identical results; therefore, only a single set of results is reported. For the noisy case, we report results for both the weak-form WLaSDI surrogate and the strong-form LaSDI surrogate.

The table lists the optimal parameter value $\hat{\mubf}$, the corresponding objective and constraint values ($f$ and $g$), and the computational cost in terms of function evaluations and runtime. All optimizations are performed on a single CPU core. For reference, we also estimate the computational cost that would be required if the optimization were performed directly using the full-order model (FOM), based on the number of function evaluations required in the surrogate-based optimization.

\begin{table}[htbp]
\centering
\renewcommand{\arraystretch}{1.1}

\begin{tabular}{lccc}
\toprule
& \textbf{Noise-free} & \textbf{WLaSDI} $\mathbf{5\%}$ \textbf{noise} & \textbf{LaSDI} $\mathbf{5\%}$ \textbf{noise} \\
%& (Differential Evolution) & (SLSQP) \\
\midrule

$\muhat$ 
& [\textbf{\num[round-mode=places,round-precision=3]{0.21459723 }},\,
\textbf{\num[round-mode=places,round-precision=3]{0.13371165}}]
&[\textbf{\num[round-mode=places,round-precision=3]{0.2119964 }},\,
\textbf{\num[round-mode=places,round-precision=3]{0.13318159}}] 
&[\textbf{\num[round-mode=places,round-precision=3]{0.20054763 }},\,
\textbf{\num[round-mode=places,round-precision=3]{0.1242434}}] \\

$f(\ubft_1, \dots, \ubft_N, \ \muhat )$
& $\num[round-mode=places,round-precision=3]{0.01936216503523197}$
& $\num[round-mode=places,round-precision=3]{0.021580364260067098}$ 
& $\num[round-mode=places,round-precision=3]{0.02065169989618149}$ \\

$g(\ubft_N, \muhat)$
& $\num[round-mode=places,round-precision=3]{0.10000002922835162}$
& $\num[round-mode=places,round-precision=3]{0.10000047157245168}$ 
& $\num[round-mode=places,round-precision=3]{0.10246862579102312}$ \\

\midrule

Func. evals 
& $564$
& $363$ 
& $239$\\

%Grad. evals \tnote{a}
%& $-$
%& $90$ \\

ROM runtime (s) 
& $\num[round-mode=places,round-precision=0]{1034.059118270874}$
& $ \num[round-mode=places,round-precision=0]{725.7901341915131}$ 
& $ \num[round-mode=places,round-precision=0]{661.0969498157501}$ \\

Est. FOM cost (s)
&  $\num[round-mode=places,round-precision=1]{5.3e7}$
& --
& --
\\

\bottomrule
\end{tabular}
\caption{\textbf{Thermal Radiative Transfer:} Differential Evolution optimization results using surrogates trained with noise-free and $5\%$ noisy data. Results for the noisy case are reported for both WLaSDI and LaSDI. WLaSDI recovers a parameter much closer to the noise-free optimal parameter than LaSDI. }
\label{tab:TRT_DiffEv}
\end{table}

\begin{figure}[htbp]
    \centering
    \includegraphics[width=0.9\textwidth]{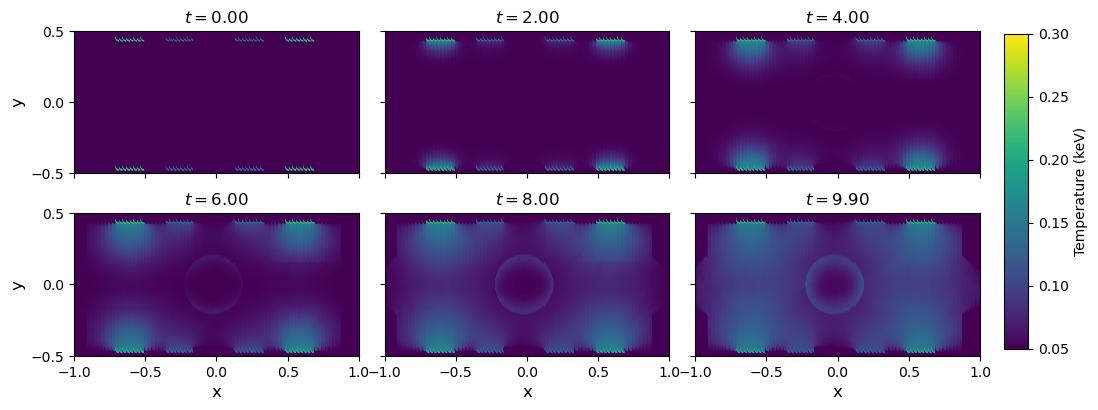}
   \caption{\textbf{Thermal Radiative Transfer:}
Time evolution of the temperature field $T(\xbf,t; \ \mubf)$ corresponding to the optimized parameter 
$\muhat = [T_{\mathrm{outer}},\,T_{\mathrm{inner}}] = [0.215,\, 0.134]$ keV, 
.}
    \label{fig:TRT_muhat_sim}
\end{figure}

\begin{figure}[htbp]
    \centering
    \begin{subfigure}{0.46\textwidth}
        \centering
        \includegraphics[width=\textwidth]{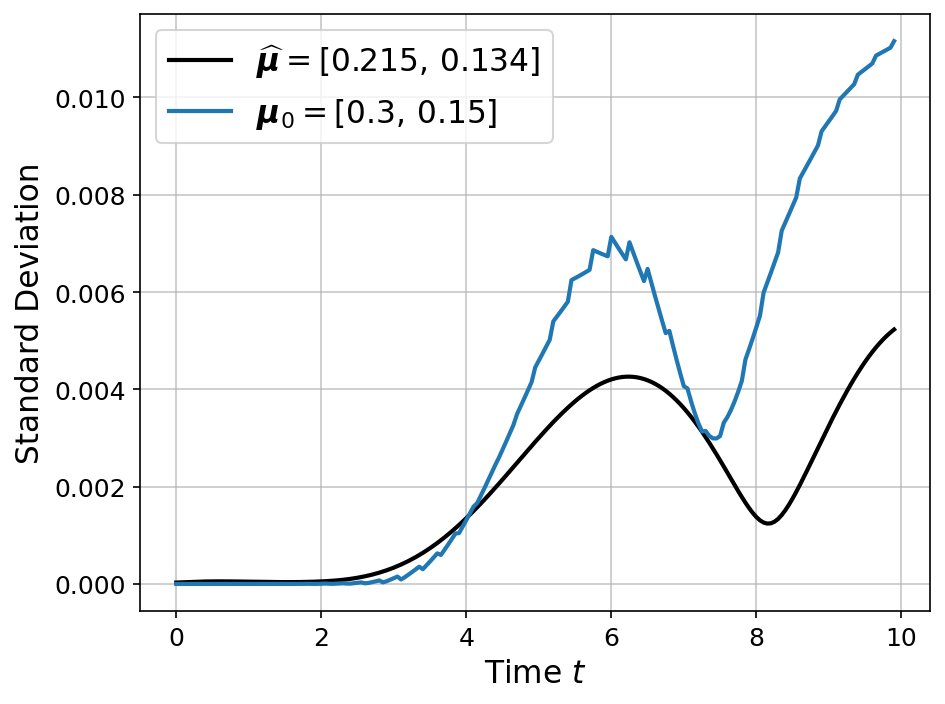}
        \caption{Temperature variation before and after optimization.}
        \label{fig:TRT_QoI_true}
    \end{subfigure}
    \hspace{0.02\textwidth}
    \begin{subfigure}{0.46\textwidth}
        \centering
        \includegraphics[width=\textwidth]{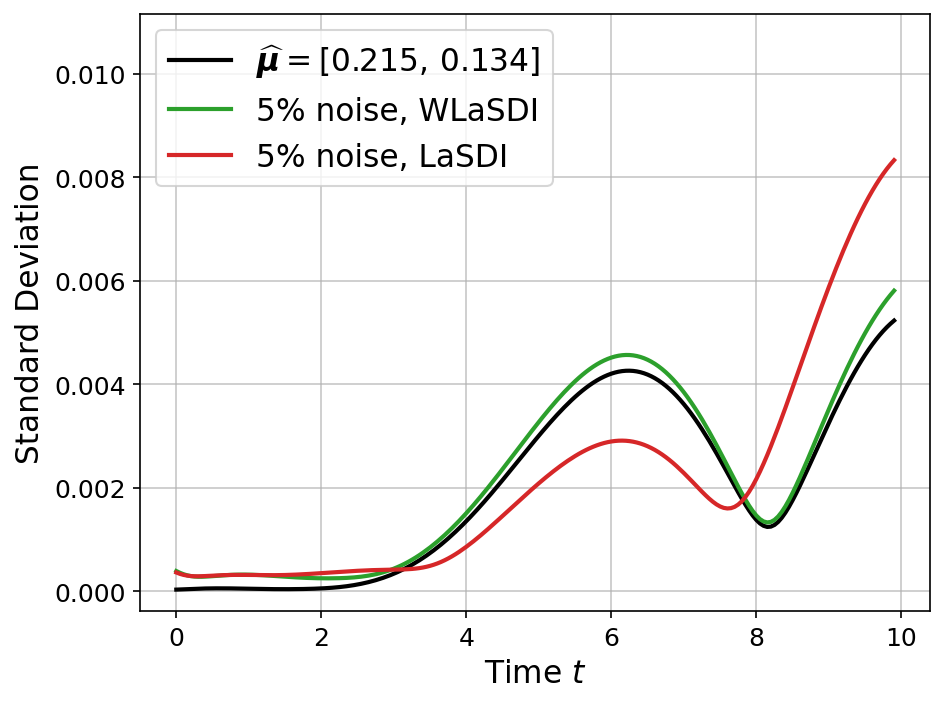}
         \caption{Comparison of the QoI by WLaSDI and LaSDI trained with $5\%$ noise.}
\label{fig:TRT_QoI_compare}
    \end{subfigure}
\caption{\textbf{Thermal Radiative Transfer:} Quantity of interest (QoI): standard deviation of the capsule surface temperature.
(a) QoI from FOM with $\mubf=[0.3,\,0.15]$ (blue) and from surrogate using noise-free global optimized parameter $\muhat=[0.215,\,0.134]$ in Eq.~\eqref{eq:TRT_muhat} (black), showing reduced temperature variation.
(b) Comparison of QoI predictions from surrogates trained with $5\%$ noisy data: noise-free reference (black), WLaSDI (green), and LaSDI (red). The WLaSDI result remains close to the noise-free solution, whereas the LaSDI prediction deviates noticeably.}
  \label{fig:TRT_QoI}
\end{figure}

Using the noise-free surrogate yields the global optimal parameter
\begin{equation}
    \muhat = [T_{\mathrm{outer}},\, T_{\mathrm{inner}}] 
= [\num[round-mode=places,round-precision=3]{0.21459723 }, \, \num[round-mode=places,round-precision=3]{0.13371165}].
\label{eq:TRT_muhat}
\end{equation}
Figure~\ref{fig:TRT_muhat_sim} illustrates the time evolution of the temperature field corresponding to the optimized parameter $\hat{\mubf}$. Figure~\ref{fig:TRT_QoI_true} compares the standard deviation of the temperature along the capsule boundary for a non-optimized parameter ($\mubf=[0.3,\,0.15]$) and for the optimized parameter $\hat{\mubf}$. The optimized configuration produces a clear reduction in temperature variation along the capsule surface.

To further assess the robustness of the weak-form, Fig.~\ref{fig:TRT_QoI_compare} compares the capsule temperature variation predicted by surrogates trained with $5\%$ noisy data. The WLaSDI surrogate (green) closely follows the noise-free reference solution (black), whereas the LaSDI surrogate (red) exhibits a noticeable deviation.

%In terms of computational cost, surrogate-based optimization yields a dramatic reduction in runtime. The optimization using the surrogate models requires only minutes on a single CPU core, whereas performing the same optimization directly with the full-order model would require an estimated runtime on the order of $10^7$ seconds. This corresponds to a speedup of approximately five orders of magnitude, demonstrating the effectiveness of the surrogate-based approach for accelerating PDE-constrained design optimization.

\paragraph{Parallel gradient-based optimization across submodels}
\label{par:TRT_parallel}

We next consider a parallel optimization strategy in which the optimization is performed independently on each of the $25$ local surrogate submodels, and the resulting local optima are compared to identify the global solution. 

Table~\ref{tab:TRT_parallel} reports the performance of the gradient-based optimizer Sequential Least Squares Programming (SLSQP) \cite{LawsonHanson1995} applied to both WLaSDI and LaSDI surrogates constructed from noise-free and $5\%$ noisy training data. Each run is initialized at the center of its corresponding subregion, and the $25$ optimizations are executed in parallel on $25$ CPU cores.
Table~\ref{tab:TRT_parallel} shows the optimal parameter $\hat{\mubf}$ across all subregions, together with the corresponding objective and constraint values ($f$ and $g$). We also report the total number of function and gradient evaluations across all subregions and the overall wall-clock runtime.

\begin{table}[htbp]
\centering
\renewcommand{\arraystretch}{1.2}
\begin{threeparttable}

\begin{tabular}{lccc}
\toprule
& \textbf{Noise-free} & \textbf{WLaSDI} $\mathbf{5\%}$ \textbf{noise} & \textbf{LaSDI} $\mathbf{5\%}$ \textbf{noise} \\
%& (Differential Evolution) &  (SLSQP) \\
\midrule
$\muhat$ 
& [\textbf{\num[round-mode=places,round-precision=3]{0.21479861}},\,
\textbf{\num[round-mode=places,round-precision=3]{0.13355837}}]
& [\textbf{\num[round-mode=places,round-precision=3]{0.21199626}},\, 
\textbf{\num[round-mode=places,round-precision=3]{0.13318034}}]
& [\textbf{\num[round-mode=places,round-precision=3]{0.2}},\, 
\textbf{\num[round-mode=places,round-precision=3]{0.12445844}}]\\

$f(\ubft_1, \dots, \ubft_N, \ \muhat )$
& $\num[round-mode=places,round-precision=3]{0.019361929473026694}$
& $\num[round-mode=places,round-precision=3]{0.021580108159595403}$ 
& $\num[round-mode=places,round-precision=3]{0.02064207385953322}$  \\

$g(\ubft_N, \muhat)$
& $\num[round-mode=places,round-precision=3]{0.0999999999935469}$ 
& $\num[round-mode=places,round-precision=3]{0.10000000000110168}$  
& $\num[round-mode=places,round-precision=3]{0.10257175736362417}$  \\

\midrule

Func. evals \tnote{a}
& $215$
& $256$  
& $230$  \\

Grad. evals \tnote{a}
& $90$
& $87$  
& $87$  \\

ROM runtime (s) \tnote{b}
& $\num[round-mode=places,round-precision=0]{58.33244157}$
& $\num[round-mode=places,round-precision=0]{23.05687499}$
& $\num[round-mode=places,round-precision=0]{32.27553654}$\\

\bottomrule
\end{tabular}

\begin{tablenotes}
\footnotesize
\item[a] Total over 25 subregions.
\item[b] Wall-clock time with parallel execution on 25 cores.
\end{tablenotes}

\end{threeparttable}
\caption{\textbf{Thermal Radiative Transfer:} Gradient-based (SLSQP) optimization results using WLaSDI and LaSDI surrogates trained with noise-free and $5\%$ noisy data. Optimization is performed independently on the $25$ submodels in parallel. In the noise-free case, the recovered parameter matches the global optimum obtained with Differential Evolution (Eq.~\eqref{eq:TRT_muhat}). With noisy training data, WLaSDI remains much closer to this optimum than LaSDI.}
\label{tab:TRT_parallel}
\end{table}

For the noise-free case, we observe that the optimal parameter $\hat{\mubf}$, along with the corresponding objective and constraint values ($f$ and $g$), is very close to the values obtained from the global optimization using Differential Evolution reported in Table~\ref{tab:TRT_DiffEv}, indicating that the method converges to the same global minimum. At the same time, the number of function and gradient evaluations, as well as the overall computational cost, is significantly reduced compared to Differential Evolution. This improvement reflects both the efficiency of the gradient-based method and the benefits of parallel execution across 25 cores.

Across the $25$ subregions, $7$ are fully infeasible in the sense that the constraint $g$ in Eq.~\eqref{eq:TRT_constraint} cannot be satisfied. As a result, the optimization in these subregions does not terminate successfully, since no feasible optimal solution exists. In contrast, for the feasible subregions, the optimization terminates successfully, identifying a local minimum that satisfies the constraint.
For the infeasible ones, we observe that the number of function and gradient evaluations is approximately $5-10$ times larger than for the feasible bins, which significantly increases the computational cost. This occurs because the optimizer continues searching for a feasible solution that does not exist. Detecting infeasible regions more efficiently is outside the scope of the present work; however, this observation motivates the use of direct optimization on the full surrogate. 

With $5\%$ noisy training data, WLaSDI and LaSDI converge to parameters consistent with those returned by Differential Evolution in Table~\ref{tab:TRT_DiffEv}. However, the WLaSDI solution remains significantly closer to the noise-free optimum $\muhat$, further highlighting the robustness achieved by the weak form.

\paragraph{Direct gradient-based optimization}
\label{par:TRT_direct}

We next present results for direct gradient-based optimization applied to the full surrogate model using a trust-region constrained optimizer \cite{ByrdHribarNocedal1999SIAMJOptim}. To assess robustness with respect to initialization, we perform $50$ optimization runs with random initializations within the admissible parameter domain, each executed on a single CPU core.

Table~\ref{tab:TRT_gradient} summarizes the results of these experiments using the WLaSDI surrogates. The reported optimal parameters correspond to the mean across runs, together with the associated $99\%$ confidence intervals. The same statistical summary is provided for the objective value $f$ and the constraint value $g$. Computational cost metrics, including the number of model evaluations\footnote{A single model evaluation corresponds to one call to the unified evaluation routine. For the derivative-free method, this call returns $f$ and $g$. For the gradient-based method, it additionally returns the gradients $\dd{f}{\mu_i}$ and $\dd{g}{\mu_i}$, all computed within the same call.} and runtime, are reported as the median together with the min–max range across the $50$ runs.

\begin{table}[htbp]
\centering
\renewcommand{\arraystretch}{1.1}
\begin{threeparttable}
\begin{tabular}{lcc}
\toprule
& \textbf{Noise-free} & \textbf{WLaSDI} $\mathbf{5\%}$ \textbf{noise} \\
%& (Differential Evolution) & (Trust-Region Constrained) \\
\midrule

$\muhat$ \tnote{a}
&  [\textbf{\num[round-mode=places,round-precision=3]{0.2149791449290836}} $\pm \num[round-mode=places,round-precision=3]{0.00009}$, \,
\textbf{\num[round-mode=places,round-precision=3]{ 0.13367075924214145}} $\pm \num[round-mode=places,round-precision=3]{0.0000667412525407988}$]
& [\textbf{\num[round-mode=places,round-precision=3]{0.21365263151677902}} $\pm \num[round-mode=places,round-precision=3]{0.0005877689502985559}$, \,
\textbf{\num[round-mode=places,round-precision=3]{0.13463823233112024}} $\pm \num[round-mode=places,round-precision=3]{0.00046507620702351067
}$] \\

$f(\ubft_1, \dots, \ubft_N, \ \muhat )$  \tnote{a}
& $\num[round-mode=places,round-precision=3]{0.019574841348619698} \pm
\num[round-mode=places,round-precision=3]{0.0000667412525407} $
& $\num[round-mode=places,round-precision=3]{0.022080110622388776} \pm
\num[round-mode=places,round-precision=3]{0.00016989118590869664} $  \\

$g(\ubft_N, \muhat)$  \tnote{a}
& $\num[round-mode=places,round-precision=3]{0.10008412893020534} \pm 
\num[round-mode=places,round-precision=3]{0.00004617769051097317}$
& $\num[round-mode=places,round-precision=3]{0.10088968972384407} \pm 
\num[round-mode=places,round-precision=3]{0.00029149365190246786}$  \\

\midrule

Func. evals \tnote{b}
& $13, 11 \!-\! 125$ 
& $12, 10 \!-\! 65$ \\

ROM runtime (s) \tnote{b}
& $\num[round-mode=places,round-precision=0]{14.332607984542847},  \num[round-mode=places,round-precision=0]{7.297398090362549} \!-\!
\num[round-mode=places,round-precision=0]{137.99986219406128}$
& $\num[round-mode=places,round-precision=0]{12.351335763931274},  \num[round-mode=places,round-precision=0]{8.681209087371826} \!-\!
\num[round-mode=places,round-precision=0]{70.23895072937012}$   \\

%\midrule
%Est. FOM cost (s) \tnote{c}
%& $\num[round-mode=places,round-precision=1]{4.5e7}$
%& $\num[round-mode=places,round-precision=1]{1.8e6}$  
%\\

%Est. FOM cost (h)
%& $\num[round-mode=places,round-precision=0]{1250}$
%& $\num[round-mode=places,round-precision=0]{500}$
%\\

\bottomrule
\end{tabular}
\begin{tablenotes}
\footnotesize
\item[a] Mean $\pm$ $99\%$ confidence interval.
\item[b] Median, min -- max range.
%\item[c] Median.
\end{tablenotes}
\end{threeparttable}
\caption{\textbf{Thermal Radiative Transfer:} Optimization results using a gradient-based trust-region constrained method applied to the WLaSDI surrogate with noise-free and $5\%$ noisy training data. Results are obtained from 50 random initializations and reported as mean $\pm$ 99\% confidence interval, with computational costs shown as min-max ranges. The method consistently converges to the global optimum.}
\label{tab:TRT_gradient}
\end{table}

As shown in Table~\ref{tab:TRT_gradient}, for both the noise-free and $5\%$ noisy cases, optimization using the WLaSDI surrogate consistently converges to the same global optimum reported in Table~\ref{tab:TRT_DiffEv}.
The number of model evaluations and runtime span a relatively wide range across runs. However, the median cost lies near the lower end of this range and remains significantly smaller than that observed for the parallel submodel optimization approach, even when the latter is executed across $25$ cores. The runtime is also substantially lower than that required by Differential Evolution. This indicates that most initializations converge efficiently, while a small number of slower cases account for the larger observed range.
Closer examination shows that, for some initializations, the optimizer temporarily stagnates before eventually converging to the global minimum. Further improvements in efficiency and robustness may be achieved through refinement of the optimizer settings; however, such investigations are beyond the scope of the present work.

\begin{figure}[htbp]
    \centering
    \includegraphics[width=0.85\textwidth]{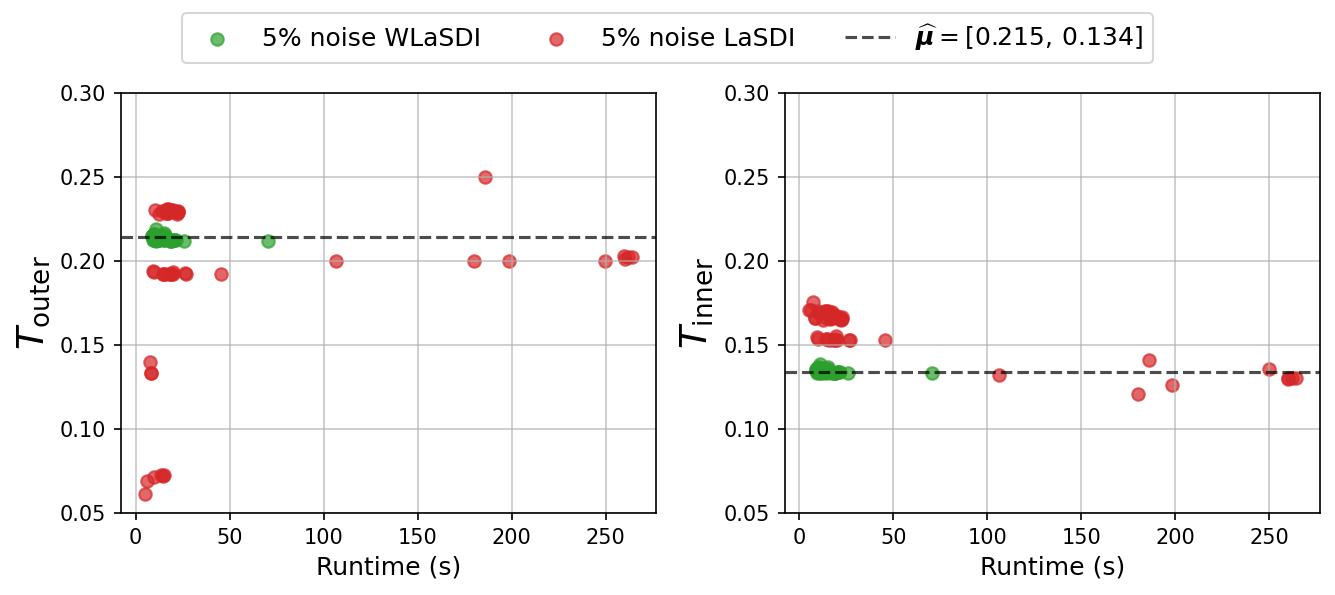}
    \caption{\textbf{Thermal Radiative Transfer:} Direct gradient-based optimization results using surrogates trained with $5\%$ noisy data. Each point corresponds to one of the $50$ runs. The dashed line indicates the noise-free global optimum $\hat{\mubf}=[0.215,\,0.134]$. WLaSDI (green) consistently recovers values close to the optimum, whereas LaSDI (red) exhibits significantly larger scatter.}
    \label{fig:TRT_gradient_compare}
\end{figure}

Finally, Fig.~\ref{fig:TRT_gradient_compare} compares the optimization outcomes obtained using WLaSDI and LaSDI surrogates trained with $5\%$ noisy data. Each point represents one of the $50$ optimization runs with a random initialization. The dashed line indicates the global optimum obtained from the noise-free surrogate in Eq.~\eqref{eq:TRT_muhat}. 
WLaSDI consistently converges to the noise-free optimum, with relatively small variability in runtime across runs. In contrast, the LaSDI exhibits significantly larger scatter in both the recovered parameters and the runtime. %These results demonstrate the robustness of using WSciML-based surrogates for accelerating gradient-based PDE-constrained optimization.

Overall, the results demonstrate that WLaSDI surrogate models provide an effective and computationally efficient framework for PDE-constrained optimization. In the TRT hohlraum design problem considered here, WLaSDI reduces the computational cost by approximately five orders of magnitude compared to optimization with the full-order model. This efficiency is further enhanced by constructing the surrogate as a collection of local submodels, which enables reliable and accurate optimization even when gradient-based methods are employed.
At the same time, the weak-form exhibits strong robustness to noisy training data, with WLaSDI consistently producing more reliable optimization results than LaSDI in the presence of noise.

The optimal hohlraum design example highlights the flexibility of WLaSDI for optimization. 
We begin with an optimization problem that minimizes the variation of the capsule surface temperature, and subsequently introduce an additional constraint to enforce physically meaningful result. More broadly, WLaSDI readily accommodates further modifications to the objective or constraints, e.g., enforcing a target heating level within a prescribed time, demonstrating its adaptability.

This underscores a key advantage of WLaSDI over surrogate approaches that approximate only a scalar QoI (e.g., radial basis function interpolation). 
By reconstructing the full solution field, WLaSDI enables the evaluation of arbitrary QoIs without retraining. This flexibility is particularly valuable in PDE-constrained optimization, where objectives and constraints may evolve over time.
In contrast, QoI-specific surrogate models must be retrained whenever a new QoI is introduced, and it can be computationally expensive, especially if the FOMs were not retained.  
By learning the full field, WLaSDI provides a reusable surrogate capable of supporting multiple downstream optimization tasks.

Taken together, these results demonstrate that WLaSDI enables significantly more efficient optimization than using full-order models,  while maintaining greater robustness and accuracy than other surrogate-based approaches.
This positions WLaSDI as a powerful and versatile tool for accelerating PDE-constrained optimization.

\section{Conclusion}
Large-scale optimization problems constrained by high-dimensional, time-dependent partial differential equations (PDEs) are computationally prohibitive, as each optimization step typically requires repeated solutions of the governing equations. The Weak-form Latent Space Dynamics Identification (WLaSDI) framework \cite{TranHeMessengerEtAl2024ComputerMethodsinAppliedMechanicsandEngineering}, a data-driven, projection-based reduced-order modeling (ROM) approach, mitigates this problem by enabling robust and efficient surrogate simulations within optimization loops.

As an extension of Latent Space Dynamics Identification (LaSDI) \cite{FriesHeChoi2022ComputerMethodsinAppliedMechanicsandEngineering}, WLaSDI projects the high-dimensional PDE onto a low-dimensional latent space, where the latent evolution is represented by a parametric ordinary differential equation (ODE). The latent dynamics are identified using Weak-form Estimation of Nonlinear Dynamics (WENDy) \cite{BortzMessengerDukic2023BullMathBiol}, which enables robust and computationally efficient simulation of high-dimensional dynamics, even with noisy training data.

In this work, we apply WLaSDI to accelerate PDE-constrained optimization. By substituting the full-order model (FOM) with a WLaSDI surrogate, governed by a low-dimensional latent ODE and coupled to the physical space through encoder–decoder mappings, the optimization no longer requires repeated high-fidelity solves. Instead, it relies entirely on the efficient solution of the latent ODE system. We present the complete formulation of this approach and derive adjoint-based sensitivities to enable efficient gradient evaluation through the latent ODE dynamics.

We demonstrate the performance of WLaSDI surrogates on several benchmark problems, including the optimal design of a hohlraum device in inertial confinement fusion governed by a thermal radiative transfer system, as well as inverse parameter recovery for the Vlasov-Poisson and Burgers' equations. We compare WLaSDI against direct optimization using full-order models (FOMs) and alternative surrogate-based approaches.
WLaSDI achieves up to five orders of magnitude speedup relative to FOM-based optimization, while maintaining high accuracy and robustness. In particular, it improves accuracy by up to two orders of magnitude compared to other surrogate methods, even when trained on noisy data.

%We further examined how modeling choices within WLaSDI influence performance in PDE-constrained optimization. Although the framework consistently remains computationally efficient and robust compared to using FOMs, configuration choices can affect accuracy and computational cost. A better understanding of these effects lies beyond the scope of the present work but would provide valuable guidance for applying WLaSDI/LaSDI in complex applications.

Despite the strong performance demonstrated in this work, the reliability of the optimization process ultimately depends on the accuracy of the WLaSDI surrogates. It is therefore crucial to understand, quantify, and control the uncertainty from the surrogates. Without a rigorous treatment of uncertainties, such surrogates cannot be reliably applied.

Another important direction for future work is the development of multi-fidelity optimization strategies that combine high-fidelity (HF) FOMs with low-fidelity (LF) WLaSDI surrogates. In this setting, WLaSDI can be used for rapid exploration, while the FOMs are invoked adaptively when higher accuracy is needed and to refine the surrogate. Systematically coupling LF and HF models within the optimization loop can improve robustness and reliability in complex problems.

\section*{Acknowledgments}
The authors would like to thank S.~Becker (University of Colorado Boulder) for insightful comments and discussions.

This work was supported by the U.S. Department of Energy, Office of Science, Office of Advanced Scientific Computing Research, through the CHaRMNET Mathematical Multifaceted Integrated Capability Center (MMICC) program.
Support was provided to the University of Colorado Boulder under Award No.~DE-SC0023346 (DMB), and to Lawrence Livermore National Laboratory (LLNL) under Award No.~DE-SC0023164 (YC).
This work was performed under the auspices of the U.S.~Department of Energy (DOE) by Lawrence Livermore National Laboratory (LLNL) under Contract No.~DE-AC52–07NA27344.
IM release: LLNL-JRNL-2017863.

%This research utilized the Blanca condominium computing resource at the University of Colorado Boulder, as well as computational resources at Lawrence Livermore National Laboratory. Blanca is jointly funded by computing users and the University of Colorado Boulder.

\bibliographystyle{elsarticle-num}   % or abbrvnat, unsrtnat, etc.
\bibliography{refs}

\appendix
\section{Interpolation of Coefficients}
\label{sec:Interpolation}

Within the LaSDI framework, latent dynamics are governed by an ODE whose coefficient matrix may depend on the physical parameter $\mubf \in \mathcal{D}$. Given training data $\Scal = \curlies{\mubf^{(k)}}_{k=1}^K$, the goal of this section is to construct an operator
$$
\Wbf(\mubf) = \boldsymbol{\Xi}\!\left(\mubf \;\middle|\; \{\Wbf^{(k)}\}_{k=1}^K, \mathcal{S}\right),
$$
that provides the coefficient matrix at arbitrary parameter values, together with its parameter sensitivities when needed.

Broadly, the approaches considered here fall into two categories:
(i) \emph{global formulations}, which identify a single parameter-independent coefficient matrix, possibly with implicit parameterization, and
(ii) \emph{explicit interpolation methods}, which construct $\Wbf(\mubf)$ directly from parameter-dependent coefficients learned at the training points.
We briefly outline these strategies below; full details are provided in \cite{BonnevilleHeTranEtAl2024}.

\subsection{Global Coefficient}
\label{sec:global}
The goal of this section is to identify a \emph{global} latent space ODE representation of the form, 
\begin{equation*}
    \dd{{\zbf}}{t}(t, \mubf) = \Wbf^T \thetabf(\zbf(t, \mubf)),
\end{equation*}
where $\Wbf$ is a parameter-independent coefficient matrix to be learned.
Here, we seek a parametric latent space ODE in which parameter variability is encoded only through the initial condition of the latent dynamics, while the governing coefficient matrix $\Wbf$ remains global and shared across all trajectories.

%Recall that in the full-order parametric PDE setting, the parameter dependence enters only through the initial condition.  By analogy, we seek a parametric latent space ODE in which parameter variability is likewise encoded solely through the initial condition of the latent dynamics, while the governing vector field remains global and shared across all trajectories.

To this end, we concatenate all training latent trajectories into a single data matrix,
$$
\Zbf := \left[\begin{array}{ccc}
         \parens{ \Zbf^{(1)}}^T &  \cdots & \parens{ \Zbf^{(K)}}^T
         \end{array}\right]^T \in \Rbb^{(N+1)K \times N_\zbf}.$$
Similarly, we construct the global feature matrix by concatenating the library evaluations $\Thetabf^{(k)}:= \Thetabf(\Zbf^{(k)})$,
$$
\Thetabf(\Zbf )  := \left[\begin{array}{ccc}
         \parens{ \Thetabf^{(1)}}^T &  \cdots & \parens{\Thetabf^{(K)}}^T
         \end{array}\right]^T \in \Rbb^{(N+1)K \times J}.
$$
Given the aggregated latent data $\Zbf$ and feature matrix $\Thetabf(\Zbf )$, the parameter independent coefficient matrix $\Wbf$ is given by  $\Wbf = \text{WENDy}(\Zbf, \Thetabf(\Zbf ) )$. Because $\Wbf$ is parameter-independent, all parameter sensitivities vanish: $\pd{\Wbf}{\mu_i} = \mathbf{0}$.

\subsection{Implicit Parameterization}
\label{sec:implicit_parameterization}
While the use of a global coefficient matrix, as described in \ref{sec:global}, is appealing, as it eliminates the need for an explicit interpolation stage and reduces computational cost, it may fail to capture inherent parameter dependence in the latent space dynamics for certain problems. To address this limitation, we introduce an implicit parameterization of the latent dynamics by augmenting the latent state with the parameter vector itself. Specifically, we define the augmentation function $\vbf: \Rbb^{N_\zbf} \times \Dcal \rightarrow  \Rbb^{(N_\zbf + N_\Dcal)}$, given by
\begin{equation*}
    \vbf(\zbf(t, \mubf), \mubf) := \left[\begin{array}{c}
       \zbf(t, \mubf) \\
         \mubf
         \end{array}\right]  \in \Rbb^{(N_\zbf + N_\Dcal)}. 
\end{equation*}
Rather than modeling the dynamics directly in terms of $\zbf$, we seek a latent ODE for the augmented variable 
$\vbf$, 
\begin{equation*}
   \dd{\vbf}{t}(t, \mubf)  
         = \Wbf^T \thetabf \parens{\vbf(t, \mubf) },
\end{equation*}
where $
\thetabf(\cdot)$ is now evaluated on the augmented latent state. The same procedure outlined in \ref{sec:global} is then applied, but with $\vbf$ replacing $\zbf$. In particular, we concatenate the augmented latent trajectories to form $\Vbf \in \Rbb^{(N+1)K \times (N_\zbf + N_{\Dcal})}$, and construct the corresponding feature matrix $\Thetabf(\Vbf) \in \Rbb^{(N+1)K \times J}$. The coefficient matrix $\Wbf \in \Rbb^{J \times (N_\zbf + N_\Dcal)}$  is subsequently identified via $\Wbf = \text{WENDy}(\Vbf, \Thetabf(\Vbf ) )$.  This formulation yields a single global model with constant coefficients, while implicitly encoding parameter dependence through the augmented state. As a result, the latent dynamics can adapt to variations in $\mubf$ without requiring parameter-dependent coefficients or post hoc interpolation. Notably, $\pd{\Wbf}{\mu_i} = \mathbf{0}$.

%As a result, the derivative of $\Wbf$ with respect to $\mubf$ is effectively zero, substantially reducing the computational cost associated with gradient evaluation

\subsection{Radial Basis Function}
\label{sec:inpterp_rbf}
We next consider explicit interpolation of parameter-dependent coefficients using radial basis functions (RBFs). Independent coefficient matrices
$\{\Wbf^{(k)}\}_{k=1}^K$ are first identified at the training parameters.
%In this section, we describe how parameter-dependent coefficient matrices are constructed via radial basis function (RBF) interpolation. This approach allows the latent dynamics to vary explicitly with the parameter through a family of coefficient matrices $\curlies{ \Wbf^{(k)}}^K_{k = 1}$.
 %$\curlies{\Wbf^{(1)}, \Wbf^{(2)}, \dots \Wbf^{(K)}}$.
To enable interpolation, we first vectorize each coefficient matrix $\Wbf^{(k)} \in \Rbb^{J \times N_\zbf}$ into a vector $\wbf^{(k)} \in \Rbb^{J N_\zbf}$. For a parameter $\mubf \in \Dcal$, the interpolated coefficient vector is expressed as
\begin{equation}
    \begin{aligned}
        \wbf(\mubf) = \sum_{k = 1}^{K} \boldsymbol{\alpha}^{(k)}  \phi \parens{\norm{\mubf - \mubf^{(k)}}_2}
    \end{aligned}
    \label{eq:RBF}
\end{equation}
where $\phi$ denotes a chosen radial basis function, such as a Gaussian or multiquadric kernel. The interpolation weights $ \curlies{\boldsymbol{\alpha}^{(k)}}_{k=1}^K \subset \Rbb^{J N_\zbf}$  are determined by enforcing exact reconstruction of the training coefficients, i.e., by solving Eq.~\eqref{eq:RBF} at the training parameters $\curlies{\mubf^{(k)}}$. Let $r^{(k)} = \norm{\mubf - \mubf^{(k)}}_2 $, it holds that
\begin{equation*}
    \begin{aligned}
        \pd{\wbf}{\mu_i}(\mubf) = \sum_{k = 1}^{K} \boldsymbol{\alpha}^{(k)} \phi' \parens{r^{(k)}} \parens{\frac{\mu_i - \mu^{(k)}_i}{r^{(k)}}} \in \Rbb^{JN_\zbf}.
    \end{aligned}
\end{equation*}
The resulting gradient vector is reshaped back into matrix form to yield $\frac{d\Wbf}{d\mu_i} \in \Rbb^{J \times N_\zbf}$.

\subsection{Convex Interpolation via Mahalanobis Distance}
\label{sec:interp_convex}
We consider a convex interpolation strategy based on Mahalanobis distance in parameter space. For a query parameter $\mubf$ and training parameters $\{\mubf^{(k)}\}_{k=1}^K$, define
\[
r^{(k)}(\mubf)
=
\bigl((\mubf-\mubf^{(k)})^{\top}\Sbf^{-1}(\mubf-\mubf^{(k)})\bigr)^{1/2},
\]
where $\Sbf$ is the empirical covariance matrix of the training parameters.
The coefficient matrix is interpolated as a convex combination of the training coefficients,
$$
\Wbf(\mubf)
=
\sum_{k=1}^{K}
\beta^{(k)}(\mubf)\,\Wbf^{(k)}, 
\qquad
\beta^{(k)}(\mubf)
=
\frac{r^{(k)}(\mubf)^{-2}}{\sum_{j=1}^{K} r^{(j)}(\mubf)^{-2}},
$$
which ensures $\sum_k \beta^{(k)}=1$.
Parameter sensitivities follow directly from differentiation of the weights,
$$
\pd{\Wbf}{\mu_i}(\mubf)
=
\sum_{k=1}^{K}
\pd{}{\mu_i} \brackets{\beta^{(k)} (\mubf)} \,\Wbf^{(k)}, 
$$
with
\begin{equation*}
    \begin{aligned}
        \pd{}{\mu_i} \brackets{ r^{(j)}(\mubf)^{-2}}
         = -2 \frac{ \parens{\Sbf^{-1} \parens{\mubf -\mubf^{(j)}}}_i}
        {r^{(j)}(\mubf)^4}.
    \end{aligned}
\end{equation*}

\subsection{Gaussian Process Regression Interpolation}
\label{sec:interp_gp}
We next consider the Gaussian Process Regression (GPR) approach for interpolating the entries of the coefficient matrix $\Wbf(\mubf)$. First introduced in GPLaSDI \cite{BonnevilleChoiGhoshEtAl2024ComputerMethodsinAppliedMechanicsandEngineering}, GPR provides a probabilistic model that yields both a mean prediction and an uncertainty estimate at each parameter value.

Each entry of $\Wbf(\mubf)$ is interpolated independently using a Gaussian process prior. Without loss of generality, let $w(\mubf)$ denote one such entry, and let $\overline{w}(\mubf)$ denote its predictive mean. The mean function is expressed as
\begin{equation*}
    \begin{aligned}
        \overline{w}(\mubf) = \sum_{k=1}^{K} \alpha^{(k)} k \parens{\mubf, \mubf^{(k)}} =  \sum_{k=1}^{K} \alpha^{(k)} \gamma \text{exp}\parens{-\frac{\norm{\mubf - \mubf^{(k)}}^2_2}{2\lambda^2}}.
    \end{aligned}
    \label{eq:GP}
\end{equation*}
Here,  $k \parens{\cdot, \cdot}$ is the covariance kernel, $\gamma $ denotes the kernel amplitude, and $\lambda  > 0$ is the characteristic length scale. The coefficients $\curlies{\alpha^{(k)}}_{k=1}^K$ are determined by solving the associated Gaussian process regression system for each of $\overline{w}^{(k)}$. 
The derivative of the predictive mean with respect to the $i$-th parameter component is then
\begin{equation*}
    \begin{aligned}
        \frac{d\overline{w}}{d \mu_i} = \sum_{k=1}^{K} \alpha^{(k)} \frac{\partial k}{\partial \mu_i} (\mubf, \mubf^{(k)})  = 
        \sum_{k=1}^{K} \alpha^{(k)} 
        \frac{\gamma}{\lambda^2}\text{exp}\parens{-\frac{\norm{\mubf - \mubf^{(k)}}^2_2}{2\lambda^2}}(\mu^{(k)}_i - \mu_i).
    \end{aligned}
\end{equation*}
The resulting partial derivatives are then reshaped and assembled across all entries to yield the full Jacobian $\pd{\Wbf}{\mu_i}(\mubf)$.

\section{Computing Derivatives for Explicit Runge--Kutta Schemes}
\label{sec:RK}
This section derives the residuals and their derivatives required for gradient-based optimization when the latent dynamics are discretized using an explicit Runge--Kutta (RK) time integration scheme. These expressions are used in Section~\ref{sec:WLaSDI_optimization} to compute parameter sensitivities and adjoint variables.

%\subsection{Runge--Kutta Time Discretization}
Recall the latent space initial value problem in Eq.~\eqref{eq:paramtric_latent_ode}
\begin{equation*}
\begin{aligned}
     \dd{\zbf}{t}(t,\mubf) &= \Wbf^T(\mubf)\thetabf\big(\zbf(t,\mubf)\big), \\
     \zbf(0,\mubf) &= \Gen\big(\gbf(\mubf)\big),
\end{aligned}
\end{equation*}
where $\zbf: [0, T] \times \Dcal \rightarrow \Rbb^{N_\zbf}$ denotes the latent state. For notational simplicity, the parametric dependence is omitted whenever it is clear from context. Let $t_n = n\dt$, with $\dt = \frac{T}{N}$ and $n = 0,1,\dots,N$.
An explicit Runge-Kutta scheme of order $s$  advances the latent state according to
\begin{equation*}
\zbf_n = \zbf_{n-1} + \dt \sum_{j=1}^s b_j \kbf_j(\zbf_{n-1},\mubf),
\qquad n = 1,2,\dots,N,
\end{equation*}
where the stage vectors $\kbf_j : \Rbb^{N_\zbf} \times \Dcal \rightarrow \Rbb^{N_\zbf}$ are defined recursively as
\begin{equation*}
\begin{aligned}
\kbf_1(\zbf,\mubf) &= \Wbf^T(\mubf)\thetabf(\zbf), \\
\kbf_j(\zbf,\mubf) &= \Wbf^T(\mubf)\thetabf\!\left(
\zbf + \dt \sum_{i=1}^{j-1} a_{ji}\ \kbf_i(\zbf,\mubf)
\right),
\qquad j=2,\dots,s.
\end{aligned}
\end{equation*}
The coefficients $\{a_{ji}, b_j\}$ correspond to the Butcher tableau of the chosen explicit Runge--Kutta method \cite{Iserles2008}.
The latent residuals associated with the Runge--Kutta discretization are enforced by the constraints
$\rbft_n = \mathbf{0}_{N_\zbf}$ for $n=0,1,\dots,N$, where
\begin{equation*}
\begin{aligned}
\rbft_0(\zbf_0,\mubf) &= \zbf_0 - \Gen\big(\gbf(\mubf)\big), \\
\rbft_n(\zbf_{n-1},\zbf_n,\mubf) &= \zbf_n - \zbf_{n-1}
- \dt \sum_{j=1}^s b_j \ \kbf_j(\zbf_{n-1},  \mubf),
\quad n=1,\dots,N.
\end{aligned}
\end{equation*}

%\subsection{Derivatives of the Residuals}
We now derive the derivatives required for sensitivity analysis and adjoint-based optimization.
The total derivative of $\rbft_0(\zbf_0, \mubf)$ with respect to $\mu_i$ is
\begin{equation*}
\begin{aligned}
\dd{\rbft_0}{\mu_i}  \big(\zbf_0, \ \mubf \big)
        &= \pd{\rbft_0}{\mu_i} \big(\zbf_0, \ \mubf \big)
        + \pd{\rbft_0}{\zbf_0} \big(\zbf_0, \ \mubf \big) 
        \ \pd{\zbf_0}{\mu_i} (\mubf)\\
        & = \parens{- \nabla \Gen(\gbf(\mubf)) 
        \ \pd{\gbf}{\mu_i}(\mubf)}
        + \parens{\Ibf_{N_\zbf} \ \pd{\zbf_0}{\mu_i} (\mubf)}.
\end{aligned}
\end{equation*}
\begin{comment}
\begin{equation*}
\dd{\rbft_0}{\mu_i}  \big(\zbf_0, \ \mubf \big)
        = \pd{\rbft_0}{\mu_i} \big(\zbf_0, \ \mubf \big)
        + \pd{\rbft_0}{\zbf_0} \big(\zbf_0, \ \mubf \big) 
        \ \pd{\zbf_0}{\mu_i} (\mubf),
\end{equation*}
with
\begin{equation*}
\begin{aligned}
 &\pd{\rbft_0}{\zbf_0}  \big(\zbf_0, \ \mubf \big) =  \Ibf_{N_\zbf}, \\
        &\pd{\rbft_0}{\mu_i} \big(\zbf_0, \ \mubf \big) = - \nabla \Gen(\gbf(\mubf)) 
        \ \pd{\gbf}{\mu_i}(\mubf).
\end{aligned}
\end{equation*}
\end{comment}
Here, $\nabla \Gen \in \Rbb^{N_\zbf \times N_\ubf}$ denotes the Jacobian of the encoder map.
For $n=1,\dots,N$, the total derivative of $\rbft_n$ satisfies
\begin{equation*}
\begin{aligned}
 \dd{\rbft_n}{\mu_i} \big(\zbf_{n-1},\zbf_n, \ \mubf\big) &= \pd{\rbft_n}{\mu_i}  \big(\zbf_{n-1},\zbf_n, \ \mubf\big)
         + \pd{\rbft_n}{\zbf_{n}}  \big(\zbf_{n-1},\zbf_n, \ \mubf\big)
         \ \pd{\zbf_{n}}{\mu_i}  (\mubf)
         + \pd{\rbft_n}{\zbf_{n-1}}  \big(\zbf_{n-1},\zbf_n, \ \mubf\big)
         \ \pd{\zbf_{n-1}}{\mu_i}  (\mubf) \\
    &= \parens{
    -\dt \sum_{j=1}^s b_j \pd{\kbf_j}{\mu_i}(\zbf_{n-1},\mubf)
    } 
    + \parens{ 
     \Ibf_{N_\zbf}  } \pd{\zbf_{n}}{\mu_i}  (\mubf)
    + \parens{
     -\Ibf_{N_\zbf}
- \dt \sum_{j=1}^s b_j \pd{\kbf_j}{\zbf}(\zbf_{n-1},\mubf)}
 \pd{\zbf_{n-1}}{\mu_i}  (\mubf).  
\end{aligned}
\end{equation*}
\begin{comment}
\begin{equation*}
\begin{aligned}
 \dd{\rbft_n}{\mu_i} \big(\zbf_{n-1},\zbf_n, \ \mubf\big) = \pd{\rbft_n}{\mu_i}  \big(\zbf_{n-1},\zbf_n, \ \mubf\big)
         + \pd{\rbft_n}{\zbf_{n-1}}  \big(\zbf_{n-1},\zbf_n, \ \mubf\big)
         \ \pd{\zbf_{n-1}}{\mu_i}  (\mubf)
         + \pd{\rbft_n}{\zbf_{n}}  \big(\zbf_{n-1},\zbf_n, \ \mubf\big)
         \ \pd{\zbf_{n}}{\mu_i}  (\mubf).
\end{aligned}
\end{equation*}
The required partial derivatives are given by
\begin{equation*}
\begin{aligned}
\pd{\rbft_n}{\mu_i}
&= -\dt \sum_{j=1}^s b_j \pd{\kbf_j}{\mu_i}(\zbf_{n-1},\mubf), \\
\pd{\rbft_n}{\zbf_n} &= \Ibf_{N_\zbf}, \\
\pd{\rbft_n}{\zbf_{n-1}}
&= -\Ibf_{N_\zbf}
- \dt \sum_{j=1}^s b_j \pd{\kbf_j}{\zbf}(\zbf_{n-1},\mubf).
\end{aligned}
\end{equation*}
\end{comment}
The derivatives $\pd{\kbf_j}{\mu_i}$ and $\pd{\kbf_j}{\zbf}$ are obtained by repeated application of the product and chain rules. In particular, these derivatives involve the Jacobian $\nabla \thetabf(\zbf)$ of the nonlinear feature map $\thetabf$ with respect to the latent state, as well as the parametric derivative $\pd{\Wbf^T}{\mu_i}(\mubf)$ of the projection matrix. The computation of $\pd{\Wbf}{\mu_i}(\mubf)$ follows from the interpolation-based construction of $\Wbf(\mubf)$ and is detailed in \ref{sec:Interpolation}.

\paragraph{Direct and Adjoint Sensitivities}

The gradient of the surrogate objective in Eq.~\eqref{eq:gradient_f_surrogate} may be computed using either direct or adjoint sensitivity analysis. For direct method, the parameter sensitivities satisfy the forward recursion
\begin{equation*}
\begin{aligned}
\pd{\zbf_0}{\mu_i}
&= -\brackets{\pd{\rbft_0}{\zbf_0}}^{-1}
\pd{\rbft_0}{\mu_i}, \\
\pd{\zbf_n}{\mu_i}
&= -\brackets{\pd{\rbft_n}{\zbf_n}}^{-1}
\left(
\pd{\rbft_n}{\mu_i}
+ \pd{\rbft_n}{\zbf_{n-1}} \pd{\zbf_{n-1}}{\mu_i}
\right),
\qquad n=1,\dots,N.
\end{aligned}
\end{equation*}
For adjoint method with adjoint variables $\ldbf_n \in \Rbb^{N_\zbf}$, the backward recursion is
\begin{equation*}
\begin{aligned}
\ldbf_N^T
&= \pd{f}{\ubf_N} \nabla \Gde(\zbf_N)
\brackets{\pd{\rbft_N}{\zbf_N}}^{-1}, \\
\ldbf_n^T
&= \left(
\pd{f}{\ubf_n} \nabla \Gde(\zbf_n)
- \ldbf_{n+1}^T \pd{\rbft_{n+1}}{\zbf_n}
\right)
\brackets{\pd{\rbft_n}{\zbf_n}}^{-1},
\qquad n=N-1,\dots,0.
\end{aligned}
\end{equation*}

\end{document}